\newif\ifHYPER\global\HYPERtrue
\newif\ifINTERNAL\global\INTERNALfalse
\documentclass[10pt,b5paper,fleqn]{article}
\usepackage[utf8]{inputenc}
\usepackage{amsmath,amsfonts,amssymb}
\usepackage[mathscr]{eucal}
\usepackage{color}
\usepackage{graphicx}
\usepackage{amsthm}
\usepackage{verbatim}
\usepackage{listings}

\ifHYPER
\definecolor{ks-green}{rgb}{0.0,0.7,0.0}
\definecolor{ks-red}{rgb}{0.7,0.0,0.0}
\definecolor{ks-blue}{rgb}{0.0,0.0,0.7}
\usepackage{hyperref}
\hypersetup{colorlinks=true,citecolor=ks-green,linkcolor=ks-red}
\fi

\numberwithin{equation}{section}
\theoremstyle{plain}

\newtheorem{proposition}[equation]{Proposition}

\theoremstyle{definition}
\newtheorem{definition}[equation]{Definition}

\newtheorem{Example}[equation]{Example}
\newtheorem{Remark}[equation]{Remark}
\newenvironment{smremark}{\emph{Remark.}}{}

\newcommand{\tsfrac}[2]{\textstyle{\frac{#1}{#2}}}

\newcommand{\tabstrut}{\rule[-0.5ex]{0pt}{3.0ex}}

\newcommand{\orth}{\bot}
\newcommand{\iprod}[2]{\langle #1, #2 \rangle}
\newcommand{\fronorm}[1]{\lVert #1 \rVert}
\newcommand{\fronormbig}[1]{\bigl\lVert #1 \bigr\rVert}
\newcommand{\norm}[1]{\lVert #1 \rVert}
\newcommand{\ledef}{\preceq}
\newcommand{\lnedef}{\prec}
\newcommand{\gedef}{\succeq}
\newcommand{\gnedef}{\succ}
\newcommand{\eigmax}{\lambda_{\text{max}}}
\newcommand{\symmat}{\mathcal{S}}
\newcommand{\sqrmat}{\mathcal{M}}
\DeclareMathOperator{\trace}{\mathrm{tr}}
\DeclareMathOperator{\diag}{\mathrm{diag}}
\DeclareMathOperator{\tomatrix}{mat}
\DeclareMathOperator{\tovector}{vec}

\def\dirsrc{.}
\def\dirfig{.}

\def\trp{^{\!\top}}
\def\minustrp{^{-\!\top}}

\def\inv{^{-1}}
\def\dd{\mathrm{d}}
\def\eps{\varepsilon}
\def\numN{\mathbb{N}}
\def\numR{\mathbb{R}}

\DeclareMathOperator{\rank}{rank}

\DeclareMathOperator{\linsp}{span}

\begin{document}
\title{Primal-dual interior-point Methods \\
  for Semidefinite Programming from an\\
  algebraic point of view, or: Using\\
  Noncommutativity for Optimization}
\author{Konrad Schrempf%
  \footnote{Contact: math@versibilitas.at (Konrad Schrempf),
    \url{https://orcid.org/0000-0001-8509-009X},
    Research Group ASiC, University of Applied Sciences Upper Austria,
    Ringstraße~43a, 4600 Wels;
    Faculty of Mathematics, University of Vienna,
    Oskar-Morgenstern-Platz~1, 1090 Wien; Austria.
    }
  \hspace{0.2em}\href{https://orcid.org/0000-0001-8509-009X}{%
  \includegraphics[height=10pt]{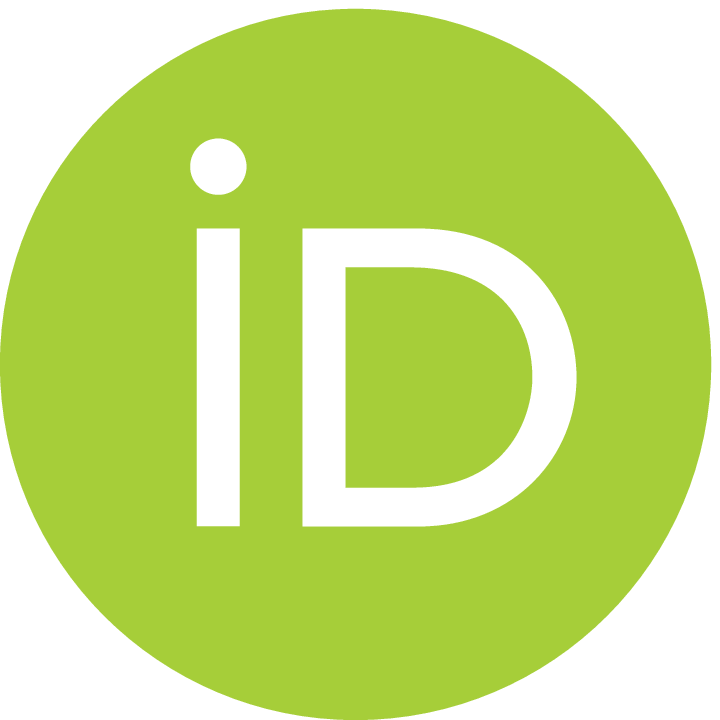}}
  }

\maketitle

\begin{abstract}
Since more than three decades, interior-point methods proved very
useful for optimization, from linear over semidefinite
to conic (and partly beyond non-con\-vex) programming;
despite the fact that already in the semidefinite case
(even when strong duality holds) ``hard'' problems are known.
We shade a light on a rather surprising restriction in the
non-commutative world (of semidefinite programming),
namely ``commutative'' paths and propose a new
family of solvers that is able to use the full richness
of ``non-commutative'' search directions:
(primal) \emph{feasible-interior-point methods}.
Beside a detailed basic discussion,
we illustrate some variants of ``non-commutative'' paths
and provide a simple implementation
for further (problem specific) investigations.
\end{abstract}

\medskip
\emph{Keywords and 2010 Mathematics Subject Classification.}
Semidefinite programming, linear programming,
interior-point methods,
primal-dual central path, sum-of-squares;
Primary 90C22, 90C51;
Secondary 90C05

\section*{Introduction}

Roughly speaking, \emph{semidefinite programming} (SDP)
is \emph{linear programming} (LP) where the (non-negative)
vectors are replaced by (positive-semidefinite) matrices.
Indeed, it is \emph{the} natural generalization of LP
by dropping the restriction to \emph{diagonal} matrices.
From a purely algebraic point of view, linear programming
is just a \emph{commutative} version of semidefinite
programming which becomes visible in particular in
(primal-dual) \emph{interior-point methods}.

\begin{smremark}
This is clearly a simplified point of view since
semidefinite programming is deeply connected with
\emph{linear matrix inequalities} (LMI's) 
\cite{Vinnikov2012b}% MR2962792 incollection
. However, this would go much too far here
and is not necessary in a first reading.
For those who are familiar with LMI's:
Is there an interpretation for the ideas
(for example \emph{algebraic barriers}
in Remark~\ref{rem:nc.algbarrier})
presented here?
\end{smremark}

\medskip
The seed of this work is one trivial observation:
Suppose that $X$ and $Z$ are \emph{invertible} (square) matrices,
say over the real numbers,
and let $I$ denote the identity matrix (all of the same size).
For a scalar $\mu > 0$ we consider the equation
\begin{displaymath}
XZ = \mu I.
\end{displaymath}
Then clearly $Z = \mu X\inv$ and since $X X\inv = X\inv X = I$ we
can conclude that $XZ = ZX$, that is, the two matrices $X$ and $Z$
\emph{commute}.
Now, for \emph{diagonal} matrices (in linear programs)
assuming commutativity is definitely not a restriction.
But what are the consequences for \emph{solving}
general semidefinite programs (by interior-point methods)?

One ``indirect'' implication is discussed in 
\cite[Section~6.1]{Todd2001a}% MR2009698 0962-4929
, namely the need for \emph{symmetrization} (of search
directions). For a more detailed discussion we recommend
\cite{Todd1999a}% MR1777451 1055-6788
. (Notice the ambiguity of the word ``symmetrization''; here
it is meant as reformulation such that the solution \emph{is}
symmetric.)
We will go a complete different way here by
``decoupling'' the matrices from the primal and the dual problem.

\medskip
As humus one can use some curiosity and 
the following remark of Helton and Putinar,
which might be a starting point for further (maybe problem specific)
development and in particular for robust and highly efficient implementations.
Once one realizes, \emph{why} no-one writes such a solver, it is usually
too late \ldots

``A lament is that all current computational semi-algebraic geometry
projects use a packaged semi-definite solver, none write their own. This
limits efficiencies for sum of squares computation.''
\cite[Section~7.1]{Helton2006c}% X0612103 arxiv

\medskip
Those who are new to semidefinite programming and
want to grow their own plant are highly
encouraged to start with 
\cite[Section~6.1]{Helton2006c}% X0612103 arxiv
\ and take the iteration scheme from
\cite[Page~6]{Alizadeh1998a}% MR1636549 1052-6234
\ as a template for an implementation.
In parallel ---for a detailed introduction, a bigger context
and further references---
we recommend to have 
\cite{Nemirovski2007a}% MR2334199 incollection
, \cite{Todd2001a}% MR2009698 0962-4929
, \cite{Vandenberghe1996a}% MR1379041 0036-1445
\ and/or 
\cite{Wolkowicz2000a}% MR1778223 book
\ at hand.
Additionally a classical book on \emph{linear programming}
like \cite{Dantzig2016a}% PUP2016 book
\ can be very useful.

\medskip
In Section~\ref{sec:nc.nla} we summarize the very basics
of (numerical) linear algebra and settle a notation which
simplifies the following discussion.
Then we recall briefly \emph{linear programming} in Section~\ref{sec:nc.lp}
before we prepare the framework for the investigation of
\emph{path-following methods} for
semidefinite programming in Section~\ref{sec:nc.sdp}
from an applied point of view (almost forgetting how they are derived).
The main contribution is the family of (primal)
\emph{feasible-interior-point methods} in Section~\ref{sec:nc.fip}
which can follow ``non-commutative'' paths (meaning that
there is no assumption on the commutation of matrices).
How to find a \emph{feasible} starting point is somewhat
rudimentary
sketched in Section~\ref{sec:nc.fim}.
The illustrations are based on a \emph{concrete problem}
discussed in Section~\ref{sec:nc.sos} and a \emph{simple implementation}
(in \textsc{Octave}) in Section~\ref{sec:nc.imp}
which is intended to serve as a starting point for
further investigations.
Section~\ref{sec:nc.stud} is an appetizer for semidefinite
programming especially for students.
That there is quite a lot of topics we can hardly touch
becomes clear from the literature which is \emph{not} cited
here. Since this work is only the very beginning we
point out some directions for further development and
(try to) mention at least starting points to facilitate
digging into different areas in Section~\ref{sec:nc.epi}.

\begin{figure}[ht]
\begin{center}
\includegraphics{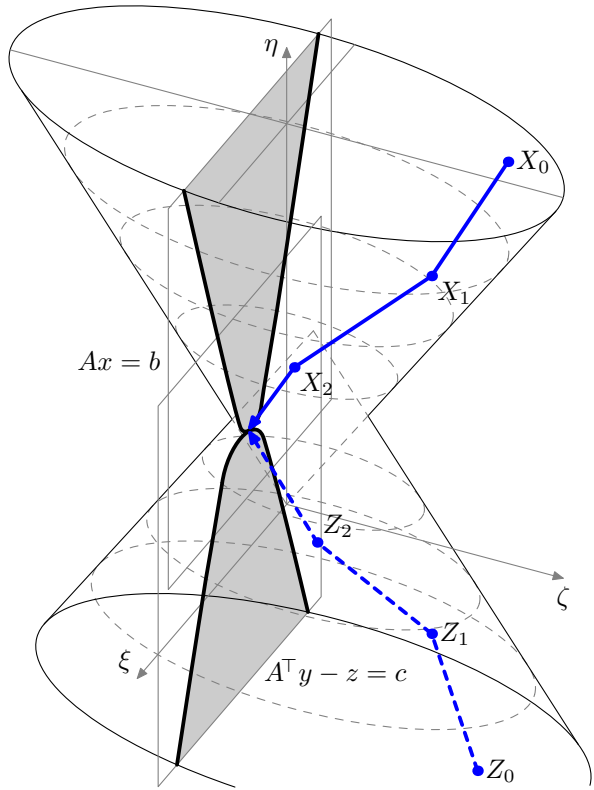}
\end{center}
\vspace{-2ex}
\caption{The ``primal'' (respectively ``dual'') cone
  of \emph{symmetric} positive semidefinite
  matrices $X\in \symmat_n$ (respectively $Z \in \symmat_n$)
  intersected with the ``hyperplane'' of
  constraints $A\tovector X = b$ (respectively $A\trp y - \tovector Z = \tovector C$).
  The primal-dual interior-path starts at $(X_0,y_0,Z_0)$
  and approaches the \emph{feasible} optimum.
  (Notice that the $X_i$'s and $Z_i$'s live in different coordinate systems
  and the illustration is with respect to the emphasis of a
  zero duality gap.)
  }
\label{fig:nc.sdp.103}
\end{figure}

\medskip
\begin{smremark}
Those who expect highly sophisticated analytical tools and
investigations might be disappointed. The only prerequisites
are literacy in \emph{linear algebra}, some experience
in \emph{numerics} (mainly from an applied point of view)
and a tiny fraction of \emph{algebra}.
For an immediate comparison of the basic idea compare
Figure~\ref{fig:nc.sdp.103} (page~\pageref{fig:nc.sdp.103})
and Figure~\ref{fig:nc.sdp.104} (page~\pageref{fig:nc.sdp.104}).
\end{smremark}

\begin{smremark}
Although we give some information on the number of iterations
in comparison with some classical solvers to stress the
---maybe surprising--- simplicity of the presented concept,
we do \emph{not} claim anything about usability and/or ``speed''
in general. The provided code is intended for
\emph{educational purposes}, in particular for
analyzing concrete optimization problems
to be able to develop more sophisticated (and optimized)
solvers, maybe in combination with existing ones.
\end{smremark}

\begin{smremark}
When we use the word ``path'' here, it is usually meant
in a general sense since we leave the ``commutative''
(analytic/continuous) \emph{central path} and need to find some
other orientation in the high-dimensional vector space
of (square) matrices to be able to get to
the minimum.
Having the \emph{polyhedron} from a linear program
and the path (from vertex to vertex) ---generated
by an instance of the Simplex Method--- \emph{on} it
in mind, that generated
by \emph{feasible-interior-point methods}
can take ``shortcuts'' \emph{through} the interior
of the polyhedron (respectively \emph{spectrahedron}
in the SDP case).
Since it is not clear a priori how to avoid
``detours'' (or oscillations)
like those shown in Figure~\ref{fig:nc.sdp.202}
(page~\pageref{fig:nc.sdp.202})
a lot of questions arise immediately.
Some of them are adressed in Section~\ref{sec:nc.epi}.
\end{smremark}

\section{(Numerical) Linear Algebra}\label{sec:nc.nla}

As a warm-up we fix the notation and recall
some basic linear algebra. Although there is nothing
really difficult, \emph{applied} (numerical)
linear algebra can be very subtle.
But we should not worry too much since we
can rely on classical functions and algorithms
from \textsc{Lapack} \cite{LAPACK2018}% ??? manual
\ which are heavily used in 
\textsc{Octave} \cite{OCTAVE2018}% ??? manual
\ and \textsc{Matlab} \cite{MATLAB2018}% ??? manual
. For a thorough discussion of \emph{applied numerical linear algebra}
we refer to \cite{Demmel1997a}% MR1463942 book
, for linear algebra in semidefinite programming to
\cite[Chapter~2]{Wolkowicz2000a}% MR1778223 book
. In Section~\ref{sec:nc.epi} we point out
more specialized literature.

Except for the notion of ``centering a matrix'' ---which has a
natural geometric interpretation (see for example Figure~\ref{fig:nc.sdp.104})
and is highly relevant from a numerical point of view---
there is nothing new here.
However, one should recall the difference between
the \emph{vector space} of (real square) matrices $\numR^{n \times n}\simeq\numR^{n^2}$
and the \emph{algebra} of matrices (with the non-commutative
multiplication) $\sqrmat_n(\numR)$.

\medskip
The set of the natural numbers is denoted by $\numN = \{ 1, 2, \ldots \}$,
that of the real numbers by $\numR$. Zero entries in matrices
are usually replaced by (lower) dots. By $I_n$ we denote the
identity matrix (of size $n$) respectively $I$ if the size is clear
from the context. We use capital letters (mainly $X$ and $Z$)
for matrices in $\sqrmat_n = \sqrmat_n(\numR)$ and use
their respective lower case letters usually to denote their ``vectorized''
version in $\numR^{n^2}$, that is,
$x = \tovector X$ (the first $n$ components in $x$ are the first
column of $X$, that from $n+1$ to $2n$ the second column, etc.),
or the other way around, $X = \tomatrix x$.
There is one exception to be consistent with the classical
notation in \emph{linear programming} in Section~\ref{sec:nc.lp},
namely $X = \diag x$ for the diagonal matrix constructed by
a vector $x \in \numR^{n \times 1}$, that is, $x_{ii}=x_i$
for $X=(x_{ij})$. The advantage is to be able to view a
\emph{linear program} (LP) as a ``commutative''
\emph{semidefinite program} (SDP)
since $\diag x \diag z = \diag z \diag x \in \sqrmat_n$.

\begin{smremark}
At a first glance it seems a little bit strange to talk about
a matrix-matrix multiplication in linear programming (which
we recall in the next section). But since our focus is with
respect to interior-point methods (mainly discussed in Section~\ref{sec:nc.sdp})
our (sequence of) vector(s) $x$ is \emph{positive}, that is,
all $n$ components of $x$ are positive and therefore in particular
invertible. In other words: There exists a vector $\tilde{x}\in \numR^{n}$
such that $\diag x \diag \tilde{x} = I_n$.
\end{smremark}

\medskip
We are mainly working with \emph{symmetric} matrices
$X = X\trp \in \symmat_n = \symmat_n(\numR) \subsetneq \sqrmat_n
\simeq \numR^{n\times n}$.
A (not necessarily symmetric) matrix $X \in \sqrmat_n$ is
called \emph{positive definite} (written as $X \gnedef 0$)
if $y\trp X y > 0$ for all $y \in \numR^n$
and \emph{positive semidefinite} (written as $X \gedef 0$)
if $y\trp X y \ge 0$ for all $y \in \numR^n$.
The \emph{trace} $\trace: \sqrmat_n \to \numR$
induces the \emph{inner product}
$\iprod{.}{.}: \sqrmat_n \times \sqrmat_n \to \numR$,
\begin{displaymath}
\iprod{X}{Z} := \trace(X\trp Z)
  = \sum_{j=1}^n \sum_{i=1}^n x_{ij} z_{ij}
  \quad\text{``$=$''}\quad (\tovector X\trp)\trp \tovector Z
\end{displaymath}
for $X = (x_{ij})$ and $Z = (z_{ij})$.
By $X \orth Z$ we denote that $X$ and $Z$ are \emph{orthogonal},
that is, $\iprod{X}{Z}=0$.

A \emph{symmetric} matrix $X=X\trp$ has
$n$ (not necessarily distinct) real \emph{eigenvalues}
$\lambda_{\min}(X) = \lambda_1 \le \lambda_2 \le \ldots \le \lambda_n
= \lambda_{\max}(X) \in \numR$.
Our main concern is the decomposition of $\sqrmat_n$ (as a vector space)
in $\ell$ subspaces $\sqrmat_n^i$ of dimension $m_i$, that is,
$\sqrmat_n = \sqrmat_n^1 \oplus \sqrmat_n^2 \oplus \ldots \oplus \sqrmat_n^\ell$
(with $m_1+m_2+\ldots+m_\ell=n^2$). As indices we will use mostly
$i \in \{ \xi, \eta, \zeta, \nu \}$ and define
the (vectorsub-)space of \emph{non-symmetric matrices} as
the orthogonal complement of the symmetric matrices,
that is,
$\sqrmat_n^\nu = (\symmat_n)^\orth
= \linsp M_\nu$
with \emph{basis} $M_\nu \in \numR^{n^2 \times m_\nu }$
(of $m_\nu=n(n-1)/2$ column vectors).
With a decomposition of $\symmat_n$ we have
\begin{displaymath}
\sqrmat_n = \underbrace{\sqrmat_n^\xi \oplus \sqrmat_n^\eta \oplus
  \sqrmat_n^\zeta}_{=\symmat_n} \oplus \sqrmat_n^\nu
\end{displaymath}
with basis $M = [M_\xi,M_\eta,M_\zeta,M_\nu]$
and $n^2 = m_\xi + m_\eta + m_\zeta + m_\nu$.
Notice that some $m_i$ can be zero
and recall that $\rank M_i = m_i$.
By abuse of notation we write
$X \in M_i = \bigl\{ M_i^{(1)},\ldots,M_i^{(m_\xi)}\bigr\}$
for $X = \tomatrix M_i^{(j)}$ for some $1 \le j \le m_\xi$.

\begin{smremark}
The decomposition of the vector space of symmetric matrices
into ``degrees of freedom'' $\sqrmat_n^\xi$,
``minimization directions'' $\sqrmat_n^\eta$ and
``constraints'' $\sqrmat_n^\zeta$ is used in particular
in Section~\ref{sec:nc.fip},
illustrated in Figure~\ref{fig:nc.sdp.104} (page~\pageref{fig:nc.sdp.104}).
It is a reminiscence of (the graph of) a \emph{concave}
function in the $x$-$y$-plane (with a possible parameter
in $z$-direction).
\end{smremark}

\medskip
For a matrix $X =(x_{ij}) \in \sqrmat_n$ we denote by
$\fronorm{X} = \sum x_{ij}^2 = \iprod{X}{X}$
the \emph{Frobenius norm}
which is just the ``classical'' vector norm
$\fronorm{X}= \norm{\tovector X}_2$.
Very often we need to ``update'' $X$,
written as $X + \alpha \Delta X$ with ``direction''
$\Delta X \in \sqrmat_n$ (``$\Delta X$'' is one symbol
and we write $\Delta x$ for ``$\Delta\tovector X$'')
and ``steplength'' $\alpha \in \numR$.
When $X$ is \emph{symmetric} and \emph{positive definite},
that is $X=X\trp \gnedef 0$ and $\Delta X \in \symmat_n$
we denote by
\begin{displaymath}
[0,\infty] \ni \alpha_{\max}^\pm = \alpha_{\max}^\pm(X,\Delta X)
  = \sup \bigl\{\alpha \in \numR \mid X \pm \alpha \Delta X \gnedef 0 \bigr\}
\end{displaymath}
the \emph{maximal} ``steplength(s)''.
We write $\alpha_{\max}$ for $\alpha_{\max}^+$.
Furthermore $X=X\trp \gnedef 0$ admits the \emph{Cholesky factorization}
$X = L L\trp$ (with lower triangular $L$).
It can be used to compute the maximal steplength
\cite[Section~2]{Alizadeh1998a}% MR1636549 1052-6234
:
\begin{equation}\label{eqn:nc.nla.steplength}
\alpha_{\max} = \bigl( \eigmax(-L\inv \Delta X L\minustrp) \bigl)\inv
\end{equation}
which is simply $\alpha_{\max} = 1/\max \diag (-X\inv \Delta X)$
for diagonal matrices $X$ and $\Delta X$.
For $\tau < 1$ we can ensure
$X + \tau \alpha_{\max} \Delta X \gnedef 0$.
However, from a numerical point of view,
we could run into troubles if $\tau = 1 - \eps$ for $\eps \ll 1$
since $\fronorm{X} \ll \fronorm{X\inv}$.
Now let $\sqrmat_n^\xi = \linsp M_n^\xi$
such that $X_\xi \orth \Delta X$
for all $X_\xi \in M_\xi$
and $\beta^\pm = \alpha_{\max}^\pm(X,X_\xi)$
assuming $0 < \beta^\pm < \infty$.
If we have the freedom to replace $X$ by some $X' = X + X_\xi$
(strictly) ``between''
$X - \beta^- X_\xi$ and
$X + \beta^+ X_\xi$
such that $\fronorm{(X')\inv} \ll \fronorm{X\inv}$,
we should do that. We will call that
``centering'' of $X$ (with respect to $X_\xi$).
We will refer to
minimizing $\fronormbig{(X+X_\xi)\inv}$ for $X_\xi \in \sqrmat_n^\xi$
as ``algebraic'' centering and
$X + \tsfrac{1}{2}(\beta^+ - \beta^-)X_\xi$
as ``geometric'' centering.

\begin{smremark}
The ``symmetrization'' $\frac{1}{2}(X+X\trp)$
of a \emph{non-symmetric} matrix $X \in \sqrmat_n$
can be interpreted as ``arithmetic'' centering
(with respect to $\sqrmat_n^\nu$).
Since usually there should not be that much confusion
between ``purely'' non-symmetric matrices and those
which have a non-trivial symmetric part, we use the
adjective ``non-symmetric'' in a sloppy way to simplify
the wording.
Otherwise one could alternatively use the adjective
``\emph{a}symmetric''.
\end{smremark}

\begin{definition}[Geometric and Algebraic Centering of a Matrix]
\label{def:nc.centering}
Let $0 \lnedef X \in \symmat_n$ and $M_\xi$ a basis of
$\sqrmat_n^\xi \subsetneq \symmat_n$
with $\dim \sqrmat_n^\xi = m_\xi \ge 1$.
If $\alpha_{\max}^+(X, X_\xi)$ and
$\alpha_{\max}^-(X, X_\xi)$ are both \emph{finite}
for all $ X_\xi \in M_\xi$ then
\begin{displaymath}
X_{\text{cen}} = X_{\text{cen}}(M_\xi)
  = \frac{1}{2} \sum_{X_\xi \in M_\xi}
  \bigr( \alpha_{\max}^+X_\xi
  - \alpha_{\max}^-X_\xi \bigr)
\end{displaymath}
(respectively $X_{\text{cen}} = X$ for $m_\xi=0$)
is called a \emph{geometric center} of $X$ with
respect to $M_\xi$.
Any $X'$ with
\begin{displaymath}
\fronormbig{(X')\inv} = \min_{X_\xi \in \sqrmat_n^\xi} \fronormbig{(X+X_\xi)\inv}
\end{displaymath}
is called an \emph{algebraic center} of $X$ with
respect to $\sqrmat_n^\xi$.
\end{definition}

\begin{Example}\label{ex:nc.geocen}
Let $n=2$ and
\begin{displaymath}
X =
\begin{bmatrix}
1 & 0.9 \\
0.9 & 1
\end{bmatrix},
\quad
\Delta X =
\begin{bmatrix}
-1 & . \\
. & -1
\end{bmatrix},
\quad
X_\xi =
\begin{bmatrix}
. & 1 \\
1 & . 
\end{bmatrix}.
\end{displaymath}
Then $\alpha_{\max} = 0.1$ and $X + \alpha_{\max} \Delta X$ is \emph{singular}.
However
\begin{displaymath}
\fronormbig{\bigl(X + (1-10^{-k}) \alpha_{\max} \Delta X\bigr)\inv}\approx 10^{k+1}
\end{displaymath}
while
$\fronormbig{(X + \alpha_{\max}\Delta X - 0.9 X_\xi)\inv}\approx 1.5713$.
Now let $k=10$ and $\alpha = (1-10^{-k})\alpha_{\max}$.
Then $\beta^+ = \alpha_{\max}^+(X,X_\xi) \approx 10^{-11}$
and $\beta^- = \alpha_{\max}^-(X,X_\xi) \approx 1.8$.
Let $\beta = \frac{1}{2}(\beta^+ - \beta^-)$.
Now the \emph{geometric} center is
$X_{\text{cen}} \approx X + \alpha \Delta X + \beta X_\xi$.
\end{Example}

\begin{Remark}
Since our main purpose is to stay away from the ``boundary''
(of singular matrices) rather than finding a precise center
we did not yet make that rigorous since we have
---for practical reasons---
an \emph{orthogonal} basis in mind.
One would expect that (geometric) centering is \emph{independent}
of a particular basis. But this needs to be shown.
\end{Remark}

Usually we are interested in a convex combination of
a matrix $X=(x_{ij})$ and an approximate center $X'$, that is,
$\mu X + (1-\mu)X'$ for $0 \le \mu \le 1$. This can
be accomplished easily by solving a \emph{linear} system
of equations (assuming $\fronorm{X\inv} \gg 1$ and $\fronorm{X} \approx 1$).
By ``$\otimes$'' we denote the Kronecker tensor product,
for example
\begin{displaymath}
X \otimes Z = 
\begin{bmatrix}
x_{1,1} Z & x_{1,2} Z & \ldots & x_{1,n} Z \\
x_{2,1} Z & x_{2,2} Z & \ldots & x_{2,n} Z \\
\vdots & \vdots & \ddots & \vdots \\
x_{n,1} Z & x_{n,2} Z & \ldots & x_{n,n} Z
\end{bmatrix}.
\end{displaymath}
For $Z \gnedef 0$ the matrix $I \otimes Z$ has full rank.
Let $\tilde{Z} = Z\inv$.
Given a basis $M_\xi$ of $\sqrmat_n^\xi$ (with $m_\xi \ge 1$)
and assuming $\fronorm{\tilde{Z}}/\fronorm{Z} \gg 1$,
the \emph{least squares solution} of
\begin{equation}\label{eqn:nc.algcen}
\begin{bmatrix}
(I \otimes \tilde{Z}) M_\xi & I \otimes Z
\end{bmatrix}
\begin{bmatrix}
\hat{z} \\ \Delta \tilde{z}
\end{bmatrix}
= (1-\mu) \tilde{z}
\end{equation}
(using for example \textsc{Lapack}/DGELS)
yields the ``approximate'' algebraic center
$Z' = Z + \tomatrix M_\xi \hat{x}$ with 
$\fronorm{(Z')\inv} \ll \fronorm{Z\inv}$ for $\mu = 0$
(and $Z' = Z$ for $\mu=1$).

\begin{Remark}\label{rem:nc.chol}
We did not investigate that in detail but it seems that
the Cholesky factorization is not as stable numerically
as computing the inverse. Therefore an ``approximate''
algebraic centering can be used if a (positive definite)
matrix is almost singular
(at least if there are some degrees of freedom).
\end{Remark}

\begin{Example}\label{ex:nc.algcen}
For Example~\ref{ex:nc.geocen} we set $Z = X + \alpha \Delta X$
and get
\begin{displaymath}
\tomatrix M_\xi \hat{z} \approx
\begin{bmatrix}
0 & -1 \\
-1 & 0
\end{bmatrix}
\quad\text{and thus}\quad
X' = 
\begin{bmatrix}
0.9 & -0.1 \\
-0.1 & 0.9
\end{bmatrix}
\end{displaymath}
with $\fronorm{(X')\inv} \approx 1.6008$.
\end{Example}

\begin{Remark}
How to solve (over- and) underdetermined linear systems
of equations is discussed in detail in
\cite[Chapter~3]{Demmel1997a}% MR1463942 book
. Given a linear system $A x = b$,
a \emph{least squares solution} $x$ satisfies
\begin{displaymath}
\norm{A x - b}_2 = \min_{\bar{x} \in \numR^n} \norm{A \bar{x} - b}_2.
\end{displaymath}
For sparse matrices (and iterative methods) we refer
to \cite{Morikuni2015a}% MR3317782 0895-4798
\ and their references.
It would be interesting to investigate the
``double approximate'' centering compared to the
linear system \eqref{eqn:nc.algcen} from before
with approximated inverse (in the context of semidefinite
programming as it is suggested in Section~\ref{sec:nc.fip}).
\end{Remark}

A set $\mathcal{K} \subseteq \sqrmat_n(\numR)$ is called \emph{convex}
if $\mu X + (1-\mu) Z$ holds for all $X,Z \in \mathcal{K}$ and $0 < \mu < 1$.
A linear equality $A x = b$ with $b \in \numR$ (and $A\in \numR^{1 \times n}$)
defines a \emph{hyperplane} in $\numR^n$, the corresponding
linear inequality $Ax \le b$ a halfspace in $\numR^n$.
The intersection of several halfspaces given by $A x \le b$ and
$b \in \numR^{m\times 1}$ for $m > 1$ is again \emph{convex}.
Therefore the intersection of $\mathcal{K}$ and $Ax\le b$ is \emph{convex}.

\section{Linear Programming}\label{sec:nc.lp}

We briefly recall the basics:
Let $1\le m<n$ be positive integers.
Given a (real) matrix $A \in \numR^{m \times n}$ (with rows $A_1,\ldots,A_m$)
and two vectors $b \in \numR^{m\times 1}$ and $c\in\numR^{n\times 1}$,
we want to find a vector $x\in \numR^{n \times 1}$ and
the \emph{minimum} of the \emph{linear objective function}
$f_0 = c\trp x$ such that the entries of $x$ are non-negative
---written as $x \ge 0$--- and
the \emph{linear constraints} $f_i = A_i x = b_i$ for
$i \in \{ 1,2,\ldots, m \}$ ---written as linear system $Ax=b$---
hold. This is called the \emph{primal problem}
\begin{equation}\label{eqn:nc.lp.primal}
\min_{ 0 \le x \in \numR^n }
  \bigl\{ c\trp x : Ax=b \bigr\}.
\end{equation}
The corresponding \emph{dual problem} is
\begin{equation}\label{eqn:nc.lp.dual}
\max_{\substack{0 \le z \in \numR^n \\ y \in \numR^m}}
  \bigl\{ b\trp y : A\trp y + z = c \bigr\}.
\end{equation}
For simplicity we assume that $\rank A = m$ and the
problems are \emph{feasible}, that is,
there exists an $x$ (respectively $y$ and $z$)
such that the minimum in $\eqref{eqn:nc.lp.primal}$ (respectively
the maximum in $\eqref{eqn:nc.lp.dual}$) is attained.

\begin{figure}
\begin{center}
\includegraphics{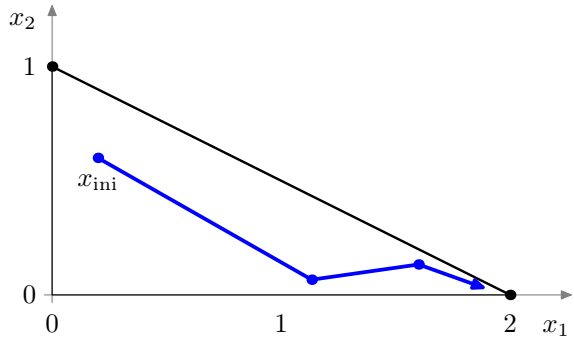}
\vspace{-1.5ex}
\caption{Tiny linear program with 
  $A = [1, 2]$, $b = [2]$, $c = [0,1]\trp$ and $x_{\text{opt}} = [2,0]\trp$.
  The ``classical'' interior path starts at $x_{\text{ini}} = [0.2,0.6]\trp$
  and approaches $x_{\text{opt}}$.
  Notice in particular that $x_{\text{ini}}$ is \emph{not} feasible.
  }
\end{center}
\label{fig:nc.sdp.101}
\end{figure}

\medskip
The ``classical'' solution is via \emph{Simplex Method}
in \emph{finitely} many steps
\cite{Dantzig2016a}% PUP2016 book
. 
Primal-dual interior-point methods move inside, that is,
$x_j >0$ and $Ax \le  b$ (respectively $z_j >0$ and $A\trp y + z \le c$),
a \emph{convex} polyhedron.
(See Figure~\ref{fig:nc.sdp.101} for an illustration
of the ``primal'' path.)
Approaching the ``boundary''
at $x_j=0$ (respectively $z_j=0$) before the optimum is reached
is punished using \emph{barrier functions}, for example
$c\trp x - \mu ( \log x_1 + \ldots + \log x_n )$
for the Lagrange multiplier $\mu > 0$.
From the necessary condition
\begin{align*}
\frac{\dd}{\dd z_k} \biggl( b\trp y + \mu \sum_{j=1}^n \log z_j \biggr)
  &= \frac{\dd}{\dd z_k} x\trp \underbrace{A\trp y}_{=c-z} + \frac{\mu}{z_k}
   = \frac{\mu}{z_k} - x_k \stackrel{!}{=} 0
\end{align*}
for an optimum we get the
\emph{complementary condition} $x_j z_j = \mu$ (for $j=1,2,\ldots,n$).
Now, starting at the triple $(x,y,z)$, we find the new \emph{search direction}
$(\Delta x, \Delta y, \Delta z)$ by solving the
\emph{linearized} system of equations
\begin{equation}\label{eqn:nc.lp.sysdir}
\begin{bmatrix}
. & A\trp & I \\
A & . & . \\
\diag{z} & . & \diag{x}
\end{bmatrix}
\begin{bmatrix}
\Delta x \\
\Delta y \\
\Delta z 
\end{bmatrix}
=
\begin{bmatrix}
c - z - A\trp y \\
b - A x \\
(\mu - x_j z_j)_{j=1}^n
\end{bmatrix}
\end{equation}
with $2n+m$ unknowns. Given a $0 < \tau < 1$,
the \emph{steplength} $\alpha\le 1$ is chosen
such that $x+\tau \alpha \Delta x > 0$ and $z + \tau \alpha \Delta z > 0$.
To make this procedure work practically,
$\mu$ must decrease (in each step) ``in the right way''.

\medskip
Since we will discuss the general case (for semidefinite
programming) in Section~\ref{sec:nc.sdp} and develop
a new approach in Section~\ref{sec:nc.fip} we only
comment a little concerning the implementation in
Section~\ref{sec:nc.imp.lp}.

Taking the classical \emph{Lagrange multiplier} approach
resulting in $x_j z_j = \mu$ one can ask what would happen
if we take $n$ \emph{independent} $\mu_j$ which could
go to zero with different ``speed''. In principle the
condition $x_j z_j = \mu > 0$ just tells us that
both components, $x_j$ and $z_j$, should be invertible.
Focusing on the primal problem and introducing ``dummy''
variables $\tilde{x}_j$ for the respective inverses of $x_j$
we need to solve the linear system
\begin{equation}\label{eqn:nc.lp.sysnew}
\begin{bmatrix}
A & . \\
c\trp M & . \\
\diag \tilde{x} & \diag x
\end{bmatrix}
\begin{bmatrix}
\Delta x \\ \Delta \tilde{x}
\end{bmatrix}
=
\begin{bmatrix}
b - A x \\ -\gamma \\ 0
\end{bmatrix}
\end{equation}
for some $\gamma > 0$. (For the general primal-dual
formulation see \eqref{eqn:nc.sdp.newdir} in the
following section.) Assuming feasibility of $x$
and a solution $\Delta x$ of \eqref{eqn:nc.lp.sysnew}
there is \emph{no} restriction for the steplength
as long as $x + \tau \alpha \Delta x$ is positive.
Thus for a \emph{feasible} (positive) starting vector
$x$ one can use a \emph{greedy} method. This is
implemented in \verb|ncminlp| in Section~\ref{sec:nc.imp.lp}.
Recall that we assume $\rank A = m$.
By $M$ we denote the orthogonal complement of $A\trp$ in
$\numR^n$. A search direction can be computed by solving
the linear system
\begin{equation}\label{eqn:nc.lp.search}
\begin{bmatrix}
c\trp M & . \\
(\diag \tilde{x}) M & \diag x
\end{bmatrix}
\begin{bmatrix}
\hat{x} \\ \Delta \tilde{x}
\end{bmatrix}
=
\begin{bmatrix}
-\gamma \\ 0
\end{bmatrix}
\end{equation}
which is much smaller than \eqref{eqn:nc.lp.sysdir}.
The search direction is $\Delta x = M \hat{x}$.
And if \eqref{eqn:nc.lp.search} is solved only
\emph{approximately} ---say by some iterative method for sparse
matrices---, it does \emph{not} have any influence
on feasibility.

A variant of the linear system \eqref{eqn:nc.lp.search} can be used
to find an initial \emph{feasible} vector $x_{\text{ini}} > 0$
by starting with $x = \tilde{x} = [1,\ldots,1]\trp$ and solving
\begin{displaymath}
\begin{bmatrix}
A & . \\
\diag \tilde{x} & \diag x
\end{bmatrix}
\begin{bmatrix}
\Delta x \\ \Delta \tilde{x}
\end{bmatrix}
=
\begin{bmatrix}
b - A x \\ 0
\end{bmatrix}
\end{displaymath}
to get a search direction $\Delta x$.
If $\alpha = 1/\lambda_n\bigl(-\diag \Delta x\, (\diag x)\inv\bigr) > 1$
we are done by using $x_{\text{ini}} = x + \Delta x$.
Otherwise we set $x: = x + \tau \alpha \Delta x$ and $\tilde{x} := x\inv$
(meaning componentwise inverse respectively the inverse as
diagonal matrix) for some $0 < \tau < 1$ and iterate.

\begin{Remark}
Table~\ref{tab:nc.lp} shows the number of iterations
for a small linear program compared with classical
solvers. Although that is rather promising, large scale
tests needs to be done. It is not yet clear how to
prove \emph{convergence} since we ``lost'' the
analytical setting. See also Section~\ref{sec:nc.imp.lp}.
For further information on Karmarkar's algorithm
\cite{Karmarkar1984a}% MR779900 0209-9683
\ in linear programming we refer to
\cite{Adler1989a}% MR1028226 0025-5610
.
\end{Remark}

\begin{Example}\label{ex:nc.lp}
Let $b = [2,7,3]\trp$, $c = [-1,-2,0,0,0]\trp$ and
\begin{displaymath}
A =
\begin{bmatrix}
-2 & 1 & 1 & . & . \\
-1 & 2 & . & 1 & . \\
1 & . & . & . & 1
\end{bmatrix}.
\end{displaymath}
Then $x_{\text{opt}} = [3,5,3,0,0]\trp$.
The following output shows that only one iteration is necessary
to find a (numerically) feasible starting vector and 3~iterations
to find the minimum (up to an error smaller than $10^{-9}$).
Notice in particular the right column with increasing norm of
the ``inverse'' of $x$ due to the last two vanishing components
in $x_{\text{opt}}$.
\end{Example}

\newpage
\begin{verbatim}
octave:1> ncminex
octave:2> x = ncminlp(A_lp, b_lp, c_lp)

NCMINLP Search STD       Version 0.99       December 2018       (C) KS
----------------------------------------------------------------------
Linear Program: n=5, m=3
      1-tau=1.00e-10, tol=2.00e-09, maxit=30
ini   alpha      min(x)        tr(c*x)             ||A*x-b|| ||x.^-1||
  1   1.80e+15    1.000e+00    -7.000000000e+00     2.55e-15   1.6e+00
cnt   alpha      tr(c*dx)      tr(c*x)             ||A*x-b|| ||x.^-1||
  0   0.00e+00   -1.000e-01    -7.000000000e+00     2.55e-15   1.6e+00
  1   3.17e+01   -3.171e+00    -1.017140874e+01     2.64e-15   5.0e+09
  2   2.83e+01   -2.829e+00    -1.300000000e+01     2.95e-15   8.7e+09
  3   3.54e+10   -4.243e-10    -1.300000000e+01     2.55e-15   7.1e+19
                               -1.299999999994142e+01

x =
   3.0000e+00
   5.0000e+00
   3.0000e+00
   5.8587e-11
   1.4143e-20

octave:3> norm(x-x_lp)
ans =    7.1754e-11
\end{verbatim}

\begin{table}
\begin{center}
\begin{tabular}{llr}
Solver & Algorithm & Iterations \\\hline\tabstrut
\textsc{SeDuMi} 1.30 \cite{Sturm1999a}% MR1778433 1055-6788
  & 0 & 17 \\
  & 1 v-corr. & 5 \\
  & 2 $xz$-corr (PC) & 3 \\\hline\tabstrut
\textsc{SDPA} 7.3.9 \cite{Yamashita2003a}% MR2019042 1055-6788
  & PC, $\tau=0.9$ & 12 \\
  & PC, $\tau=0.99999$ & 9 \\\hline\tabstrut
\texttt{ncminlp} (Section~\ref{sec:nc.imp.lp})
  & $\tau=0.99999$ & $1+5$ \\
  & $\tau=1-10^{-10}$ & $1+3$ \\\hline
\end{tabular}
\end{center}
\vspace{-2ex}
\caption{The number of iterations (that for
  finding a feasible starting vector $+$ that for minimization)
  for the linear program
  from Example~\ref{ex:nc.lp}.
  ``PC'' stands for \emph{predictor-corrector}
  which implies more expensive steps.
  For \textsc{SDPA} the sparse format with 5~diagonal
  blocks of size $1\times 1$ is used, for \textsc{SeDuMi}
  one can specify that the problem is linear.
  For both solvers their respective standard settings are used.
  One must be careful with a direct comparison due
  to the large number of possible tuning parameters.
  There are almost no parameters in \texttt{ncminlp}.
  }
\label{tab:nc.lp}
\end{table}

\clearpage

\section{Semidefinite Programming}\label{sec:nc.sdp}

Semidefinite programming is \emph{the} natural generalization
of linear programming by going over from (the \emph{cone} of)
non-negative vectors to the cone of \emph{positive semidefinite}
symmetric matrices.
The \emph{primal problem} (in standard form) is
\begin{equation}\label{eqn:nc.sdp.primal}
\min_{0 \ledef X \in \symmat_n} \Bigl\{ \iprod{C}{X} :
  \iprod{A_1}{X} = b_1, \iprod{A_2}{X} = b_2, \ldots, \iprod{A_m}{X} = b_m \Bigr\}.
\end{equation}
The corresponding \emph{dual problem} is
\begin{equation}\label{eqn:nc.sdp.dual}
\max_{\substack{0 \ledef Z \in \symmat_n,\\\atop y \in \numR^m }}
  \Bigl\{ b\trp y : y_1 A_1 + y_2 A_2 + \ldots + y_m A_m + Z = C \Big\}.
\end{equation}
Using diagonal matrices $X = \diag x$, $C = \diag c$, $Z = \diag z$
and (by abuse of notation) $A_k = \diag A_k$, a linear program
can directly be formulated as \emph{semidefinite program} (SDP).
For the connection with \emph{linear matrix inequalities} (LMI's)
we refer to 
\cite{Vinnikov2012b}% MR2962792 incollection
, for \emph{matrix inequalities} (MI's) to 
\cite{Camino2006a}% MR2219142 1052-6234
.

\begin{smremark}
The inner product $\iprod{C}{X}$ is also often denoted by
$C \bullet X$ (for \emph{symmetric} matrices $C$ and $X$).
\end{smremark}

\medskip
For convenience we state at least the definition of the
\emph{central path} and recall \emph{weak duality}.
Notice that \emph{strong duality}, that is,
$\iprod{X}{Z} = 0$, implies $XZ=0$
\cite{Todd2001a}% MR2009698 0962-4929
. Why? Does strong duality also imply $\rank X + \rank Z = n$?
For an overview about the differences between linear and
semidefinite programming see in particular
\cite[Section~3.5]{Alizadeh1995a}% MR1315703 1052-6234
. For details we refer to 
\cite{Wei2010a}% MR2718693 0025-5610
.

\begin{definition}[Central Path 
\protect{\cite[Section~3.1]{Todd1998a}% MR1636553 1052-6234
}, 
\protect{\cite[Section~10.2]{Wolkowicz2000a}% MR1778223 book
}]\label{def:nc.centralpath}
The set of solutions to the equations
\begin{align*}
\iprod{A_i}{X} &=  b_i \quad\text{for all $i \in \{1,2,\ldots, m\}$}, \\
\sum_{i=1}^{m} y_i A_i + Z &= C \quad\text{and}\\
XZ &= \mu I
\end{align*}
for all $\mu > 0$ such that $X \gnedef 0$ and $Z \gnedef 0$
is called \emph{central path}.
\end{definition}

\begin{proposition}[Weak Duality
\protect{\cite[Proposition~2.1]{Todd2001a}% MR2009698 0962-4929
}]\label{pro:nc.wdual}
If $X$ is feasible in \eqref{eqn:nc.sdp.primal} and
$(y,Z)$ is feasible in \eqref{eqn:nc.sdp.dual}, then
\begin{displaymath}
\iprod{C}{X} - b\trp y = \iprod{X}{Z} \ge 0.
\end{displaymath}
\end{proposition}

\medskip
Starting from the triple $(X,y,Z)$, a new search direction
$(\Delta X, \Delta y, \Delta Z)$
can be computed by solving a linear system of equations
similar to~\eqref{eqn:nc.sdp.sysdir}
with $2n^2 +m$ unknowns:
\begin{equation}\label{eqn:nc.sdp.sysdir}
\begin{bmatrix}
. & A\trp & I_{n^2} \\
A & . & . \\
I_n \otimes Z & . & I_n \otimes X
\end{bmatrix}
\begin{bmatrix}
\Delta x \\
\Delta y \\
\Delta z 
\end{bmatrix}
=
\begin{bmatrix}
c - z - A\trp y \\
b - A x \\
\tovector(\mu I_n - XZ)
\end{bmatrix}.
\end{equation}
Here we are only interested in the principal idea
of finding a (matrix-valued) search direction.
For detailed discussions and references we refer to
\cite{Alizadeh1998a}% MR1636549 1052-6234
\ and \cite{Todd1998a}% MR1636553 1052-6234
, with respect to (many subtle technical details for)
an implementation to
\cite{Mehrotra1992a}% MR1186163 1052-6234
. There is a whole zoo of search directions
discussed in 
\cite{Todd1999a}% MR1777451 1055-6788
.

\medskip
The last matrix equation of \eqref{eqn:nc.sdp.sysdir} comes from the
\emph{linearization} of the \emph{centering condition} $XZ = \mu I$
(for $\mu > 0$) 
\cite{Alizadeh1998a}% MR1636549 1052-6234
:
\begin{align*}
(X + \Delta X)(Z + \Delta Z) 
  &= XZ + X\Delta Z + \Delta X Z + \Delta X \Delta Z \stackrel{!}{=} \mu I.
\end{align*}
This condition implies that $Z$ is a scalar multiple of $X\inv$ and
therefore $XZ = ZX$, that is, $X$ and $Z$ \emph{commute}.
Although this is a classical assumption and a ``natural''
generalization of the complementary condition $x_j z_j = \mu$
from (primal-dual in\-te\-rior-point methods for) linear programming
we have not found much discussions about its implications
except for the need for symmetrization in 
\cite{Todd1999a,Todd2001a}% MR2009698 0962-4929
.
Even if \emph{strong duality} holds, that is,
the \emph{duality gap} is zero, 
\begin{displaymath}
\iprod{C}{X_{\text{opt}}} - b\trp y_{\text{opt}} = \iprod{X}{Z_{\text{opt}}} = 0,
\end{displaymath}
it is far from clear what happens in \emph{practical}
situations (where the optimum is never reached \emph{exactly}).
We conjecture that the origin of \emph{hardness} of SDP programs in
the sense of \cite{Wei2010a}% MR2718693 0025-5610
\ is the restriction to ``commutative'' (central) paths.
While in principle this would not be any problem
for \emph{linear programming} (and diagonal matrices)
one can use a much weaker assumption, namely that of
\emph{invertibility} (of the components), leading to
a new class of interior-point methods.
For a short discussion in the context of LP see
Section~\ref{sec:nc.lp}.

\medskip
When we go over from the \emph{analytic} formulation (in particular
of some \emph{barrier} functions 
\cite[Section~3]{Nemirovski2007a}% MR2334199 incollection
) to the \emph{discrete} (algebraic) ---not necessarily linear---
formulation, some information is ``lost''.
What we use in fact in interior-point methods
(for convex problems)
is nothing more than the possibility to find ``good''
\emph{search directions}
for minimization (respectively maximization) ---often by solving
linear systems of equations---
and the \emph{maximal steplength}
for a given search direction such that we do not
leave the ``interior''.
For semidefinite programming, both steps are rather simple
(at least if we ignore fast convergence issues),
the tricky part is numerics (and the implementation).

\medskip
Now we assume that $X$ is \emph{positive definite},
in particular \emph{invertible}, that is,
there exists an invertible $\tilde{X}$ such that
$X\tilde{X} = \tilde{X}X = I$
\cite[Remark~4.13]{Schrempf2018c}% X180905425 arxiv
.
Since we are looking for an update $\Delta X$ such that
$X +\Delta X$ is still invertible, we get the condition
\begin{align*}
(X + \Delta X)(\tilde{X} + \Delta \tilde{X})
  &= X \tilde{X} + X\Delta\tilde{X} + \Delta X \tilde{X} + \Delta X \Delta \tilde{X}
  \stackrel{!}{=} I
\end{align*}
which we linearize (assuming ``small'' $\Delta X$ and $\Delta \tilde{X}$)
and reformulate (assuming only $\Delta X \tilde{X} = \tilde{X} \Delta X$) as
\begin{equation}\label{eqn:nc.sdp.invertibility}
X \Delta\tilde{X} + \tilde{X} \Delta X = I - X\tilde{X} = 0.
\end{equation}
Thus ---compared to \eqref{eqn:nc.sdp.sysdir}---
we can find also ``non-commutative'' primal-dual search directions by solving
the following linear system of equations
(with $4n^2+m$ unknowns):
\begin{equation}\label{eqn:nc.sdp.newdir}
\begin{bmatrix}
. & A\trp & I_{n^2} & . & . \\
A & . & . & . & . \\
z\trp & . & x \trp & . & . \\
I_n \otimes \tilde{X} & . & . & I_n \otimes X & . \\
. & . & I_n \otimes \tilde{Z} & . & I_n \otimes Z \\
\end{bmatrix}
\begin{bmatrix}
\Delta x \\ \Delta y \\ \Delta z \\ \Delta \tilde{x} \\ \Delta \tilde{z}
\end{bmatrix}
=
\begin{bmatrix}
c - z - A\trp y \\
b - A x \\
\mu - \iprod{X}{Z} \\
0^{n^2 \times 1} \\
0^{n^2 \times 1} 
\end{bmatrix}.
\end{equation}
Recall that $x=\tovector X$, etc. Its system matrix has $3n^2+m+1$ rows,
that is, the system is \emph{underdetermined}
(we always assume $\rank A = m$).
The third block row reads
$\iprod{Z}{\Delta X} + \iprod{X}{\Delta Z} = \mu - \iprod{X}{Z}$.
If one is familiar with the ``classical'' issues of
primal-dual interior-point methods like symmetrization,
step-length computation, update of $\mu$, etc.~then an
implementation is ---at least for moderate sized problems---
straight forward.

\begin{Remark}
An ad hoc implementation applied to the problem
from Section~\ref{sec:nc.sos} showed that the primal path
``converges'' much faster than the dual path.
We have not investigated that in detail (for different
problem classes) since working with the \emph{primal} problem
only has the advantage of a much smaller linear system of equations
(to compute a search direction).
It would be interesting to see what happens if one uses alternating
``primal'' and ``dual'' steps,
or how they can be used to accelerate ``classical''
interior-point methods.
\end{Remark}

\begin{Remark}\label{rem:nc.algbarrier}
It is a funny coincidence that the
approach we tried 
to ensure invertibility of an (unknown) matrix $X$,
namely introducing another matrix $\tilde{X}$ and add
the equations $X \tilde{X} = I$,
as an alternative to the classical $\det X \neq 0$ 
seems to have a drawback when a linearization (like we can use here)
is not appropriate (and Groebner bases are needed) because
the number of unknowns doubles
\cite[Remark~4.13]{Schrempf2018c}% X180905425 arxiv
.
The important observation is the origin of
$XZ=\mu I$ from the classical barrier $-\log \det X$
(respectively $\log \det Z$)
\cite{Nemirovski2007a}% MR2334199 incollection
, 
\cite[Section~2]{Todd1998a}% MR1636553 1052-6234
. So it is rather natural to carry that over directly to
our discrete (iterative) setting as ``algebraic''
barrier $X \tilde{X} = I$ (respectively $Z \tilde{Z} = I$).
Recall that $\det X > 0$ does \emph{not} imply
positive eigenvalues $\lambda_i > 0$.
\end{Remark}

\newpage
\section{Feasible-interior-point Methods}\label{sec:nc.fip}

\begin{figure}
\begin{center}
\includegraphics{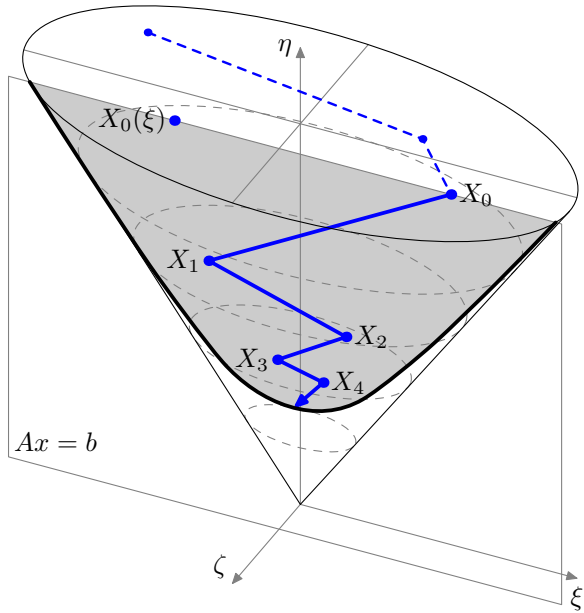}
\end{center}
\vspace{-2ex}
\caption{The cone of \emph{symmetric} positive semidefinite 
  matrices $X = X\trp \in \symmat_n \subsetneq \sqrmat_n(\numR)$
  intersected with the ``hyperplane'' of constraints $Ax=b$
  (from an SDP) for $x = \tovector X$.
  The primal \emph{feasible-interior-points}
  $(X_i)$ for $i=1,2,\ldots$~form a path towards the minimum.
  Each matrix $X_i$ satisfies $A x_i=b$,
  and can be ``moved'' along the $\xi$-axis without changing the
  value of the (linear) objective function $\iprod{C}{X_i}$,
  that is, $\frac{\dd}{\dd\, \xi} \iprod{C}{X_i(\xi)}=0$.
  The dashed (blue) line illustrates the ``initial iterates''
  to get a \emph{feasible} starting matrix $X_0$
  (discussed in Section~\ref{sec:nc.fim}).
  }
\label{fig:nc.sdp.104}
\end{figure}

Now we start (somewhat naively) from scratch
with the basic linear algebra from Section~\ref{sec:nc.nla}
in mind and recall the relevant setup for the
semidefinite program \eqref{eqn:nc.sdp.primal},
\begin{displaymath}
\min_{0 \ledef X \in \symmat_n} \Bigl\{ \iprod{C}{X} :
  \iprod{A_1}{X} = b_1, \iprod{A_2}{X} = b_2, \ldots, \iprod{A_m}{X} = b_m \Bigr\}.
\end{displaymath}
Given $m+1$ square matrices $C,A_i \in \sqrmat_n(\numR) = \numR^{n \times n}$
and the vector $b \in \numR^{m \times n}$ we want to find the
\emph{minimum} of $\iprod{C}{X} = \trace(C\trp \! X)$ such that
$X$ is \emph{symmetric} and \emph{positive semidefinite}
and the $m$ constraints $\iprod{A_i}{X} = b_i$ are fulfilled.
We write $x = \tovector X$ and
$A\trp = [\tovector A_1, \tovector A_2, \ldots, \tovector A_m]$
and use Example~\ref{ex:nc.sdp} (from the following section)
for illustration:
\begin{displaymath}
A =
\begin{bmatrix}
. & 1 & . & 1 & . & . & . & . & . \\
. & . & 1 & . & 1 & . & 1 & . & . \\
. & . & . & . & . & 1 & . & 1 & . \\
. & . & . & . & . & . & . & . & 1
\end{bmatrix},
\quad
b = \begin{bmatrix}
0 \\ \tabstrut\frac{13}{4} \\ \tabstrut\frac{15}{4} \\ 1
\end{bmatrix}
\quad\text{and}\quad
C =
\begin{bmatrix}
1 & . & . \\
. & 0 & . \\
. & . & 0 
\end{bmatrix}.
\end{displaymath}

The vector space of $n \times n$ matrices $\sqrmat_n = \sqrmat_n(\numR)$
decomposes (almost) naturally in
\emph{degrees of freedom} $\sqrmat_n^\xi$,
(pure) \emph{minimization directions} $\sqrmat_n^\eta$,
\emph{constraint directions} $\sqrmat_n^\zeta$ and
\emph{non-symmetric directions} $\sqrmat_n^\nu$,
namely
\begin{displaymath}
\sqrmat_n = \underbrace{\sqrmat_n^\xi \oplus \sqrmat_n^\eta \oplus
  \sqrmat_n^\zeta}_{=\symmat_n} \oplus \sqrmat_n^\nu
\end{displaymath}
with $m_\nu = n(n-1)/2$ and
$m_\xi + m_\eta + m + m_\nu = n^2$ and
bases $M^\xi$, $M^\eta$, $M^\zeta$ and $M^\nu$ respectively.

\begin{Remark}
However it is not so clear, how to choose $M^\eta$.
In some sense $m_\eta = \rank M^\eta$ should be \emph{maximal}
for $(M^\eta)\trp$ in row echelon form.
But on the other hand, if there are minimization
directions which ``correlate'' with other directions
it might be helpful to have more degrees of freedom
for numerical stability (at least in the beginning).
Notice that in the family of examples from
\cite{Jansson2007a}% MR2377260 0036-1429
\ $m_\xi = 0$ (see also Table~\ref{tab:nc.jck}).
\end{Remark}

Figure~\ref{fig:nc.sdp.104} shows the cone of
symmetric positive definite matrices $X = X\trp \gedef 0$
intersected with the ``hyperplane'' $Ax = b$ for $x = \tovector X$.
For a \emph{feasible}
matrix $X_k = X_k\trp \gnedef 0$,
a (symmetric) direction $\Delta X_k \in \sqrmat_n^\xi \oplus \sqrmat_n^\eta$
and a scalar $\alpha_k \in \numR$,
the matrix $X_{k+1} := X_k + \alpha_k \Delta X_k$ is \emph{again} feasible,
that is
\begin{displaymath}
A(x_k + \alpha_k \Delta x_k) = b.
\end{displaymath}
If $\Delta X_k$ is \emph{non-symmetric}, that is,
$\Delta X_k \in \sqrmat_n^\xi \oplus \sqrmat_n^\eta \oplus \sqrmat_n^\nu$,
a \emph{symmetric} direction is given by
$\Delta X_k := \frac{1}{2}\bigr(\Delta X_k + (\Delta X_k)\trp\bigl)$.
Before we discuss how to find such search directions $\Delta X_k$
(given a feasible positive definite $X_k$),
we have a look on how to compute the \emph{maximal} steplength $\alpha_k^{\text{max}}$
such that $X_k + \alpha_k^{\text{max}} \Delta X_k \gedef 0$.
(Without loss of generality we assume that $\Delta X_k$ is such that
$\alpha_k^{\text{max}} < \infty$.)

\medskip
To simplify notation we drop the iteration index in the following
and write $X = X_k$.
Since $X$ is positive definite it admits a \emph{Cholesky factorization}
$X = L L\trp$. The \emph{maximal steplength}
is \eqref{eqn:nc.nla.steplength},
$\alpha_{\max} = \bigl( \eigmax(-L\inv \Delta X L\minustrp) \bigl)\inv$
\cite[Section~2]{Alizadeh1998a}% MR1636549 1052-6234
. So for $0 < \tau < 1$ we have
\begin{displaymath}
X + \tau \alpha_{\max} \Delta X \gnedef 0.
\end{displaymath}
Now we take a (not necessarily orthonormal) basis $M^\xi$
(respectively $M^\eta$ and $M^\nu$)
for $\sqrmat_n^\xi$ (respectively $\sqrmat_n^\eta$ and $\sqrmat_n^\nu$)
and write $M^{\xi,\eta}$ (respectively $M^{\xi,\eta,\nu}$)
for $[M^\xi,M^\eta]$ (respectively $[M^\xi,M^\eta,M^\nu]$)
or simply just $M$ (for $M^{\xi,\eta}$ and $M^{\xi,\eta,\nu}$).

\begin{smremark}
If the basis $M^\xi$ is \emph{not} orthogonal one
might be a little bit more careful with respect to
\emph{geometric} centering, see Definition~\ref{def:nc.centering}
and the following remark.
\end{smremark}

\begin{smremark}
Steplength scaling can be quite subtle in practice.
For LP's, tpically $\tau=0.99995$ is used
which is too aggressive for ``predictor-corrector''
methods 
\cite[Section~4.4]{Todd1998a}% MR1636553 1052-6234
\ therefore $\tau=0.98$ (and adaptive steplength scaling)
is used.
Since here we restrict ourself to \emph{feasible} paths,
steplength scaling seems to be ---\emph{except} close to
the ``boundary'' of singular matrices--- not a difficult
issue.
\end{smremark}

\medskip
From Section~\ref{sec:nc.sdp} we recall the \emph{linearized}
condition
\eqref{eqn:nc.sdp.invertibility} for the invertibility
of $X + \Delta X$, namely
\begin{displaymath}
X \Delta\tilde{X} + \tilde{X} \Delta X = 0
\end{displaymath}
for given $X$ and $\tilde{X} = X\inv$.
Since we start with a \emph{feasible} $X$ and restrict the search
directions to $\Delta X \in \linsp M$, this invertibility condition
now reads
\begin{equation}\label{eqn:nc.fip.invertibility}
X \Delta\tilde{X} + \tilde{X} \underbrace{M \hat{x}}_{=\Delta X} = 0
\end{equation}
for $\hat{x} \in \numR^{m_\xi+m_\eta}$ or $\hat{x} \in \numR^{m_\xi+m_\eta+m_\nu}$.
Additionally, since we want to \emph{minimize},
we need a condition of the form
$\iprod{C}{\Delta X} = \iprod{C}{M\hat{x}} < 0$.
Thus, a \emph{minimization search direction} can be computed as
a \emph{least squares solution} of the (underdetermined)
\emph{linear} system of equations
\begin{equation}\label{eqn:nc.fip.search}
\begin{bmatrix}
(I_n \otimes \tilde{X}) M & I_n \otimes X \\
(\tovector C)\trp M & .
\end{bmatrix}
\begin{bmatrix}
\hat{x} \\ \Delta \tilde{x}
\end{bmatrix}
=
\begin{bmatrix}
0 \\ -\gamma
\end{bmatrix}
\end{equation}
for $\gamma > 0$. Now we can use \eqref{eqn:nc.nla.steplength}
to compute the steplength $\alpha_{\max}$ and iterate,
at least in principal.
Notice that the search direction computed in \eqref{eqn:nc.fip.search}
is \emph{independent} of $\gamma > 0$ (modulo scaling).

\medskip
\begin{smremark}
Although the path $(X_i)$ we would get ---at least if
we ignore the numerical issues---
is in the ``hyperplane'' $Ax=b$ and ``positive definite'',
we need to ensure that we indeed approach the minimum.
That this is not clear can be seen in Figure~\ref{fig:nc.sdp.202}
(cyan line): There is a first huge step, but then it seems
that the path ``converges'' to some $X'$ such that $\iprod{C}{X'}\approx 1.005 > 1$,
that is, the minimum is not reached.
One has to be very careful here because it might be that
---using multiprecision arithmetics--- just a very high number
of iterations is necessary to get arbitrary close to the minimum.
At the end this does not really matter since we want to
reduce the number of iterations as much as possible.
\end{smremark}

\begin{Remark}
For a case where centering is not needed (or even
must not be used) ---and how to adapt the implementation
from Section~\ref{sec:nc.imp}---
see (the end of) Section~\ref{sec:nc.stud}.
We have not yet investigated that in detail.
\end{Remark}

\newpage
There are roughly three different variants of
(primal) \emph{feasible-interior-point} methods, namely:

\subsection*{All Directions (STD)}

Use \emph{all} ``feasible'' search directions (including the non-symmetric),
symmetrize and center mainly for ``staying away from
the boundary'' by using a centering parameter $0.65 \le \mu \le 0.85$
and $X_{\text{new}} := \mu X + (1-\mu) X_{\text{cen}}$.
This corresponds to \verb|typ=1| in line~7
in \verb|ncminsdp| (Section~\ref{sec:nc.imp.sdp})
for \emph{algebraic} centering.

\begin{figure}
\begin{center}
\includegraphics{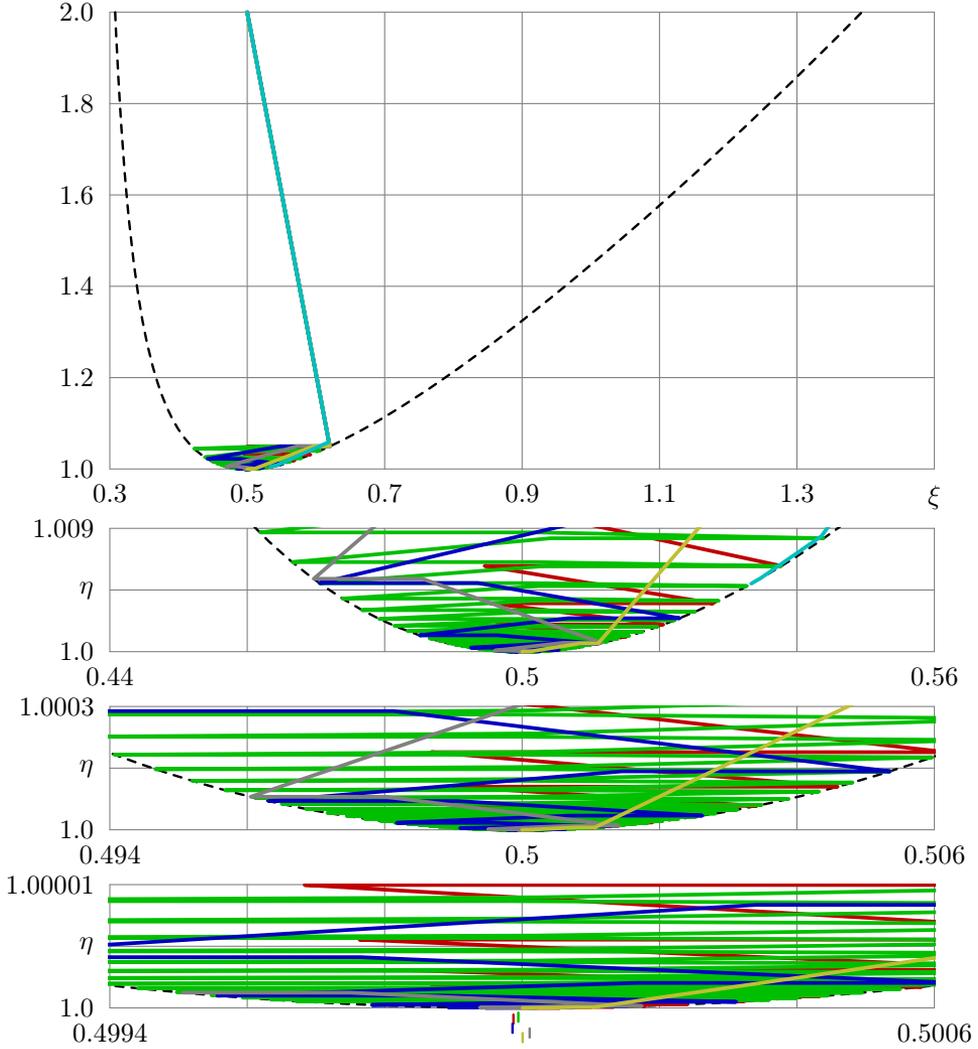}
\end{center}
\vspace{-2ex}
\caption{The paths for starting matrices $X$ with
  $\xi=0.5$, $\eta = \iprod{C}{X} = 2.0$ and
  centering parameters
  $\mu=0.0$ (red), $\mu=0.2$ (green), $\mu=0.4$ (blue),
  $\mu=0.6$ (grey), $\mu=0.8$ (yellow) and $\mu=1.0$ (cyan)
  using \emph{all} (feasible) search directions.
  The lower subfigures show a zoom around the minimum
  with $(10\times,30\times)$,
  $(100\times,900\times)$ and $(1000\times,27000\times)$
  respectively.
  In the latter, the $\xi$ values of the endpoints for
  $\mu=0.0$ (27~iterations), $\mu=0.2$ (290), $\mu=0.4$ (22), $\mu=0.6$ (9)
  and $\mu=0.8$ (6) are marked. For $\mu=0.2$ only the
  subpath of the iterations $1,2,5,8,\ldots$ is shown.
  For $\mu=1.0$ the steplength scaling has to be decreased
  (from $\tau=0.999999$) to $\tau=0.99$ to work numerically
  up to 9~iterations.
}
\label{fig:nc.sdp.202}
\end{figure}

\medskip
Some paths for different \emph{feasible} initial matrices
are shown in Figure~\ref{fig:nc.sdp.201},
some for different centering parameters in
Figure~\ref{fig:nc.sdp.202}.
For a typical problem run see (the end of)
Section~\ref{sec:nc.sos}.

Instead of the ``approximated'' algebraic centering,
``exact'' \emph{geometric centering} could be used.
It seems that the optimal centering parameter (at least for
algebraic centering) is around $\mu \approx 0.75$.
But so far it is neither clear if it is (within some range)
problem dependent nor how it can be interpreted.
A dynamical tuning (depending on the relative error)
might be possible.

\subsection*{Symmetric Directions (SYM)}

Use only \emph{symmetric} search directions $M^{\xi,\eta}$
and center mainly for ``staying close to the center''
by using a centering parameter $0 \le \mu \le 0.25$.
This corresponds to \verb|typ=2| in line~7
in \verb|ncminsdp| (Section~\ref{sec:nc.imp.sdp})
for \emph{algebraic} centering.

\begin{figure}
\begin{center}
\includegraphics{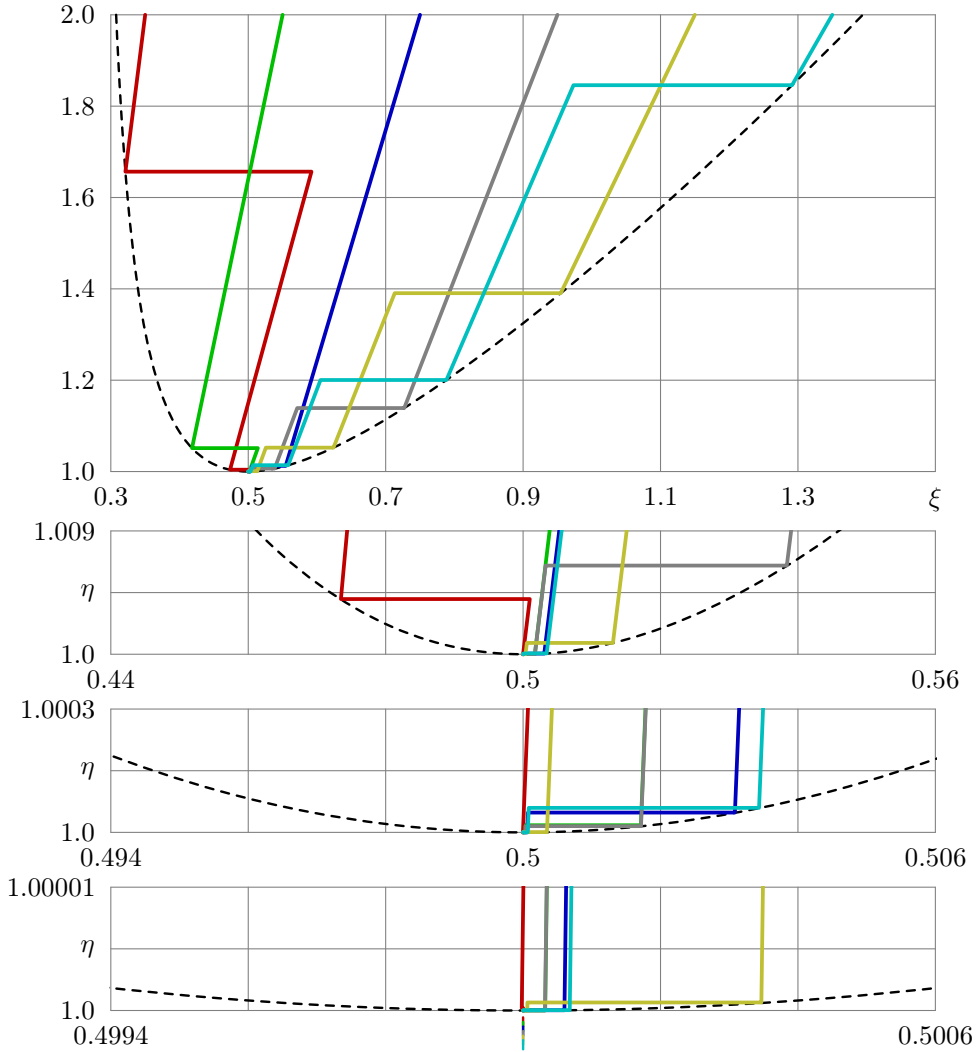}
\end{center}
\vspace{-2ex}
\caption{The paths for starting matrices $X$ with
  $\xi=0.35$ (red), $\xi=0.55$ (green), $\xi=0.75$ (blue),
  $\xi=0.95$ (grey), $\xi=1.15$ (yellow) and $\xi=1.35$ (cyan),
  each with $\eta = \iprod{C}{X} = 2.0$
  and centering parameter $\mu=0.2$
  using \emph{symmetric} (feasible) search directions only.
  The lower subfigures show a zoom around the minimum
  with $(10\times,30\times)$,
  $(100\times,900\times)$ and $(1000\times,27000\times)$
  respectively.
  In the latter, the $\xi$ values (after 4~to 7~iterations)
  of the endpoints are marked.
}
\label{fig:nc.sdp.203}
\end{figure}

\medskip
Some paths for different \emph{feasible} initial matrices are
shown in Figure~\ref{fig:nc.sdp.203}.
For a problem run (in the context of the relaxation
of a combinatorial problem) \emph{without} centering
see (the end of) Section~\ref{sec:nc.stud}.
For a comparison of the number of iterations
for the example from 
\cite[Section~6]{Jansson2007a}% MR2377260 0036-1429
\ see Table~\ref{tab:nc.jck}
(page~\pageref{tab:nc.jck}).

Instead of the ``approximated'' algebraic centering,
``exact'' \emph{geometric centering} could be used.
The latter could be used to estimate a ``good''
centering parameter $\mu$ for the former.
On the other hand, the former might be necessary
to get away from the boundary for an ``aggressive''
steplength scaling parameter $\tau < 1$ for the latter
(see Remark~\ref{rem:nc.chol}).

\subsection*{Constructed Directions (ACE or GCE)}

If degrees of freedom are available, that is,
$m_\xi \ge 1$, search directions can be constructed
by using two
\emph{different} centered matrices $X_1$ and $X_2$
with $\iprod{C}{X_2} < \iprod{C}{X_1}$
and setting $\Delta X := X_2 - X_1$.

\begin{figure}
\begin{center}
\includegraphics{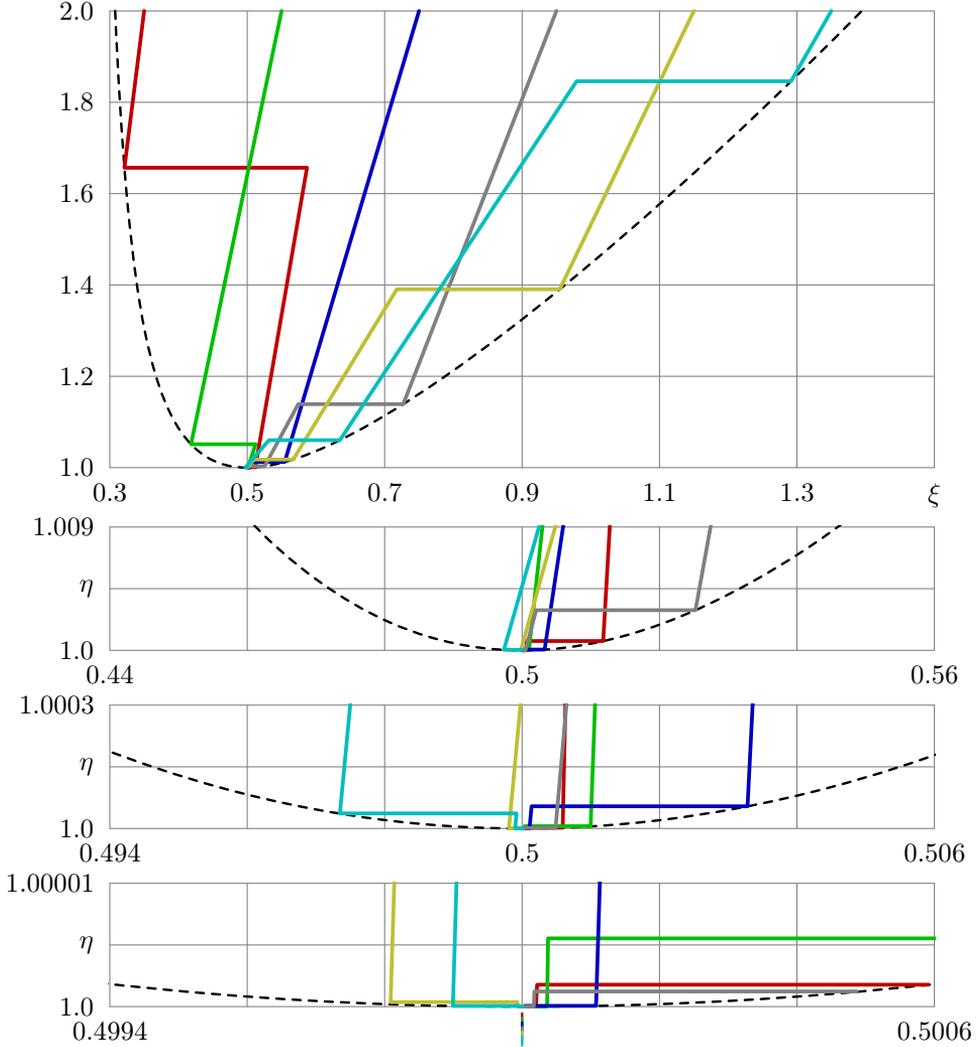}
\end{center}
\vspace{-2ex}
\caption{
  The paths for starting matrices $X$ with
  $\xi=0.35$ (red), $\xi=0.55$ (green), $\xi=0.75$ (blue),
  $\xi=0.95$ (grey), $\xi=1.15$ (yellow) and $\xi=1.35$ (cyan),
  each with $\eta = \iprod{C}{X} = 2.0$
  and centering parameter $\mu=0.215$
  using (feasible) search directions constructed by
  ``approximated'' \emph{algebraic centering}.
  The lower subfigures show a zoom around the minimum
  with $(10\times,30\times)$,
  $(100\times,900\times)$ and $(1000\times,27000\times)$
  respectively.
  In the latter, the $\xi$ values (after 4~to 6~iterations)
  of the endpoints are marked.
  Notice that in the first step \emph{symmetric} search
  directions are used, therefore the first edge of the paths
  agree with the corresponding in Figure~\ref{fig:nc.sdp.203}
  (and~\ref{fig:nc.sdp.205}).
}
\label{fig:nc.sdp.204}
\end{figure}

\medskip
For some paths using ``approximate'' \emph{algebraic}
centering see Figure~\ref{fig:nc.sdp.204}, for
some using ``exact'' \emph{geometric} centering
Figure~\ref{fig:nc.sdp.205}. An investigation for
$m_\xi \ge 2$ is open.
For a typical problem run (and some further comments)
see Section~\ref{sec:nc.sos}.
Notice however, that the usability in the case of
several minimization directions, that is, $m_\eta \ge 2$,
is not yet clear.

\begin{figure}
\begin{center}
\includegraphics{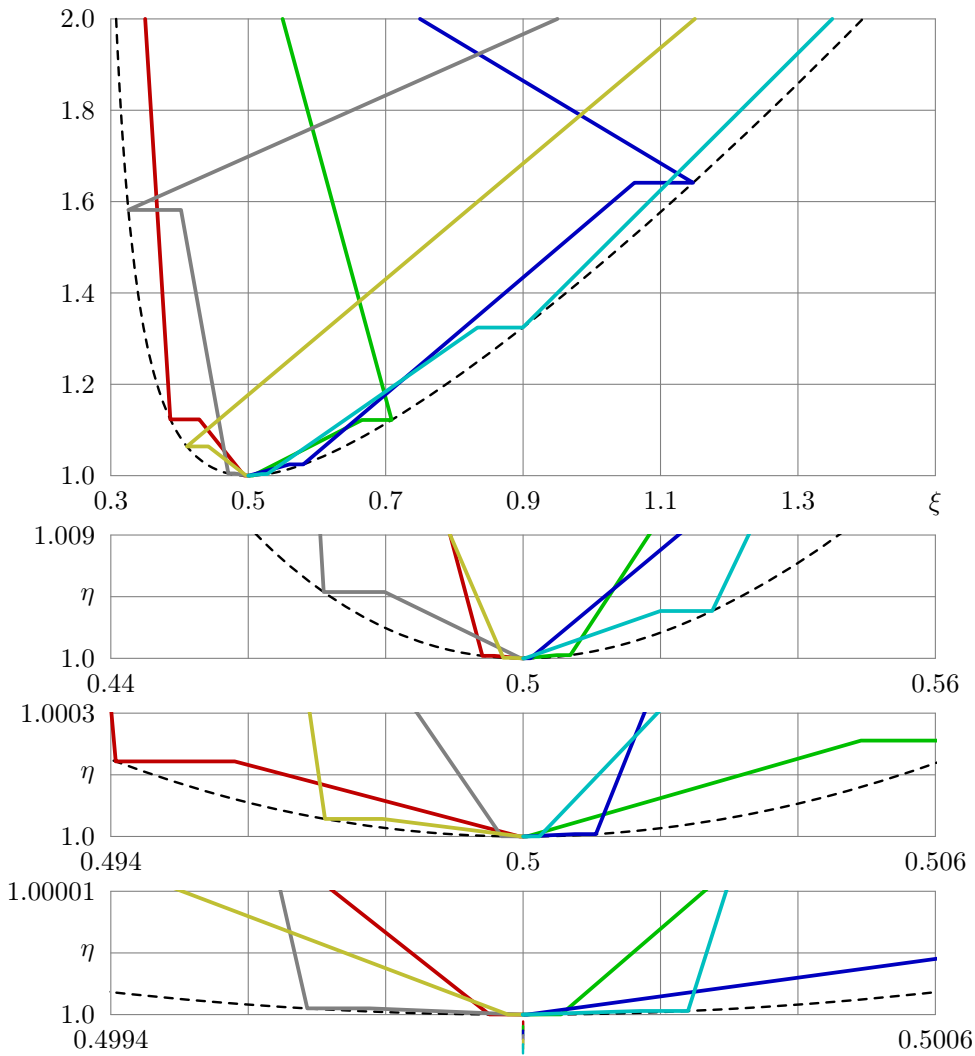}
\end{center}
\vspace{-2ex}
\caption{The paths for starting matrices $X$ with
  $\xi=0.35$ (red), $\xi=0.55$ (green), $\xi=0.75$ (blue),
  $\xi=0.95$ (grey), $\xi=1.15$ (yellow) and $\xi=1.35$ (cyan),
  each with $\eta = \iprod{C}{X} = 2.0$
  and centering parameter $\mu=0.765$ using
  \emph{all} (feasible) search directions.
  The dashed line marks the border where one
  eigenvalue becomes zero.
  The lower subfigures show a zoom around the minimum
  with $(10\times,30\times)$,
  $(100\times,900\times)$ and $(1000\times,27000\times)$
  respectively.
  In the latter, the $\xi$ values (after 4~to 5~iterations)
  of the endpoints are marked.
  Notice that one iteration corresponds to \emph{two}
  edges in the path: a ``minimization'' edge and
  a (horizontal) ``centering'' edge.
}
\label{fig:nc.sdp.201}
\end{figure}

\begin{figure}
\begin{center}
\includegraphics{\dirfig/sdp.205}
\end{center}
\vspace{-2ex}
\caption{The paths for starting matrices $X$ with
  $\xi=0.35$ (red), $\xi=0.55$ (green), $\xi=0.75$ (blue),
  $\xi=0.95$ (grey), $\xi=1.15$ (yellow) and $\xi=1.35$ (cyan),
  each with $\eta = \iprod{C}{X} = 2.0$
  and centering parameter $\mu=0.0$
  using (feasible) search directions constructed by \emph{geometric centering}.
  The lower subfigures show a zoom around the minimum
  with $(10\times,30\times)$,
  $(100\times,900\times)$ and $(1000\times,27000\times)$
  respectively.
  In the latter, the $\xi$ values (after 4~to 5~iterations)
  of the endpoints are marked.
  Notice that in the first step \emph{symmetric} search
  directions are used, therefore the first edge of the
  paths agree with the corresponding in Figure~\ref{fig:nc.sdp.203}
  (and~\ref{fig:nc.sdp.204}).
}
\label{fig:nc.sdp.205}
\end{figure}

\clearpage
\section{Finding a feasible initial matrix}\label{sec:nc.fim}

\begin{figure}
\begin{center}
\includegraphics{\dirfig/sdp.301}
\end{center}
\vspace{-2ex}
\caption{The graph of the maximal steplength $\alpha$
  from equation \eqref{eqn:nc.nla.steplength}
  of $\gamma$ along the direction $\Delta x$ as least squares
  solution of $A \Delta x = b - A \tovector(\gamma I)$
  for the example from Section~\ref{sec:nc.sos} (red)
  and the example from \cite[Section~6]{Jansson2007a}% MR2377260 0036-1429
  \ with
  $\eps=\delta=1$ (green),
  $\eps=\delta=0.01$ (blue),
  $\eps=10^{-3}$ and $\delta=-10^{-3}$ (yellow),
  $\eps=\delta=10^{-4}$ (grey) and
  $\eps=-10^{-3}$ and $\delta=10^{-3}$ (cyan).
  Notice that the latter is primal \emph{infeasible}
  and there is no $\gamma >0$ such that $\alpha(\gamma) > 1$.
  The cyan curve hides the blue, the yellow and the grey one.
  The dashed curves are scaled around $\alpha(\gamma)=1$ to
  make the shape visible.
}
\label{fig:nc.sdp.301}
\end{figure}

Even in ``classical'' interior-point methods (for semidefinite programming)
the role of ``good'' \emph{initial iterates} should not be underestimated
\cite[Section~4]{Toh1999a}% MR1778429 1055-6788
. For \textsc{SDPA} a ``sufficiently large'' multiple of the identity matrix
should be used \cite{SDPA2008}% ??? manual
, the default parameter is $\lambda^*=100$.
Since ``$\lambda$'' is heavily used for eigenvalues,
we use ``$\gamma$'' here, for example $X_0 = \gamma I$.
Unfortunately we have not found that much beyond short remarks in the literature.
For linear programming one could start with 
\cite[Section~4]{Adler1989a}% MR1028226 0025-5610
.
Since in our case finding a \emph{feasible} initial iterate is crucial,
we tried several approaches and leave some comments for further investigation.

Actually it turned out that this task seems to be highly non-trivial.
Several of our ad hoc approaches ---using variants of the techniques from
the minimization in the previous section--- failed in some cases.
So far only one survived and is presented in detail in Section~\ref{sec:nc.imp.ini}.
The variant for linear programs in Section~\ref{sec:nc.imp.lp} is just
a simplified version. If this also fails in general, a clever
use of semidefinite programming in the context of ``norm-minimization''
might work,
maybe in combination with Ger\v sgorin circles
\cite{Varga2004a}% MR2093409 book
.

\medskip
Still, starting with a least squares solution to $A x =b$ and using
repeatedly a spectral shift $\gamma I$ and a correction
$A \Delta x = b - A\bigl(x+\tovector(\gamma I)\bigr)$
seems to work for practical problems.
For small parameters $\eps$ and $\delta$ in the example
\cite[Section~6]{Jansson2007a}% MR2377260 0036-1429
\ one can increase the shift $\gamma$ dynamically.

From a different point of view one can ask how the
``projection'' of the central ray $\gamma I$ onto the
constraints $A x = b$ looks like. Not directly but in
terms of the \emph{maximal steplength} $\alpha(\gamma)$
as a function of $\gamma$, see Figure~\ref{fig:nc.sdp.301}.
If we find any $\gamma$ such that $\alpha(\gamma) > 1$ we
are done. Assuming the existence, starting with a small
$\gamma$, it is not difficult to find it with only a few
iterations (increasing $\gamma$ if the gradient is positive
and decreasing $\gamma$ if it is negative respectively
continuing with the old iterate and increase $\gamma$
less ``aggressive'').

\medskip
The search for an initial iterate is done in two steps:
Firstly, given $X$ and setting $\tilde{X} = X\inv$,
we get a direction ``towards'' $Ax=b$ as a least
squares solution to the linear system of equations
\begin{equation}\label{eqn:nc.fim.feasible}
\begin{bmatrix}
A & . \\
I_n \otimes \tilde{X} & I_n \otimes X 
\end{bmatrix}
\begin{bmatrix}
\Delta x \\ \Delta \tilde{x}
\end{bmatrix}
=
\begin{bmatrix}
b-A x \\ 0
\end{bmatrix}.
\end{equation}
Notice that a scaling of the first block row might
be necessary depending on the norm of $X\inv$.
Since $\Delta X$ is not necessarily symmetric,
we set $\Delta X := \frac{1}{2}\bigl(\Delta X + (\Delta X)\trp \bigr)$.
(In the implementation it is always done for numerical reasons.)
Now we can calculate $\alpha = \alpha_{\max}(X,\Delta X)$
and update $X := X + \min\{ 1, \tau\alpha \} \Delta X$
for some $\tau < 1$.
Secondly, for numerical stability we use
either \emph{algebraic centering} \eqref{eqn:nc.algcen}
if $m_\xi \ge 1$ or a modified version by computing
a least squares solution $\Delta X = \tomatrix(M \hat{x})$ for
$M = [ M_\eta, M_\nu ]$
to the linear system of equations
\begin{equation}\label{eqn:nc.fim.residuum}
\begin{bmatrix}
I_n \otimes \tilde{X} M & I_n \otimes X 
\end{bmatrix}
\begin{bmatrix}
\hat x \\ \Delta \tilde{x}
\end{bmatrix}
=
\begin{bmatrix}
 (1-\mu) \tilde{x}
\end{bmatrix}
\end{equation}
for $0 \le \mu \le 1$.
In both cases we update
$X := X + \frac{1}{2}\bigl( \Delta X + (\Delta X)\trp \bigr)$.
Notice that $\tilde{X} = X\inv$ needs to be updated in between.

\medskip
Using $M= M_\eta$ only in the last step did \emph{not} work
for the example from 
\cite{Jansson2007a}% MR2377260 0036-1429
. Some numbers for the initial iterates (for the presented approach)
are in Table~\ref{tab:nc.jck}.
For the problem from the next section see also
Table~\ref{tab:nc.sdp}. Notice that in this case
there is a degree of freedom for ``approximated'' algebraic
centering.

\begin{table}
\begin{center}
\begin{tabular}{llrrrr}
  & & \textsc{SeDuMi} 1.30 & \textsc{SDPA} 7.3.9 
    & \texttt{ncminsdp} & Section~\ref{sec:nc.imp.sdp} \\
$\epsilon$ & $\delta$
  & \cite{Sturm1999a}% MR1778433 1055-6788
  & \cite{SDPA2008}% ??? manual
  & standard & symmetric  \\\hline\tabstrut
$10^{-2}$ & $10^{-2}$  &  9  & 19 & $4+4$ & $4+5$ \\
$10^{-4}$ & $10^{-4}$  & 16  & 18 & $6+4$ & $6+7$ \\
$10^{-6}$ & $10^{-6}$  & 23  & 16 & $6+4$ & $6+7$ \\
$10^{-8}$ & $10^{-8}$  & (28) & 15 & $9+\text{failed}$ & $9+(6)$ \\\hline 
\end{tabular}
\end{center}
\vspace{-2ex}
\caption{The number of iterations (that for finding a feasible
  starting matrix $+$ that for minimization)
  for the semidefinite program from 
  \cite[Section~6]{Jansson2007a}% MR2377260 0036-1429
  . For \textsc{SDPA} the standard settings did not work for
  small parameters and therefore $\lambda^*=10^{6}$
  was used in all cases. 
  Iterations in parentheses mean that
  the result was not very accurate.
  }
\label{tab:nc.jck}
\end{table}

\clearpage
\section{Sum of Squares}\label{sec:nc.sos}

\begin{figure}
\begin{center}
\includegraphics[scale=1]{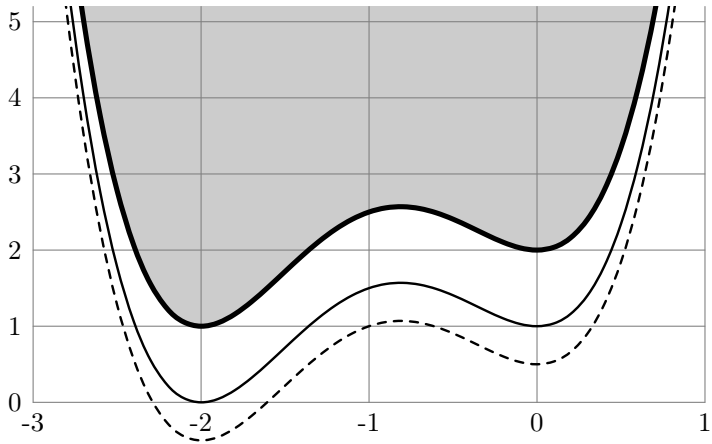}
\end{center}
\vspace{-2ex}
\caption{Graph of the polynomial
  $p_\eta(x) = \eta + \frac{13}{4} x^2 + \frac{15}{4} x^3 + x^4$ for
  $\eta = 2$ (thick line), $\eta=1$ (thin line) and $\eta=1/2$
  (dashed thin line). The polynomials $p_2$ and $p_1$
  are \emph{non-negative} and can be written as 
  \emph{sum-of-squares} (SOS), for example
  $p_2(x) = \bigl(\sqrt{2} - \frac{\sqrt{2}}{8} x^2\bigr)^2
         +\bigl(\frac{\sqrt{15}}{2} x + \frac{\sqrt{15}}{4} x^2 \bigr)^2
         +\bigl(\frac{\sqrt{2}}{8} x^2 \bigr)^2$.
  }
\label{fig:nc.poly}
\end{figure}

Suppose, that we want to find the \emph{global} minimum of the
polynomial $p(x) = 2 + \frac{13}{4} x^2 + \frac{15}{4} x^3 + x^4$
(its graph is the thick line in Figure~\ref{fig:nc.poly}).
What is easy (in this case) using high-school
mathematics will turn out to be rather involved using \emph{sum-of-squares}
and \emph{semidefinite programming} 
\cite{Lasserre2000a,Parrilo2001a,Helton2006c}% MR1814045 1052-6234
, although the main idea is so simple: Minimize $\eta$ such that
$p_\eta(x) = \eta + \frac{13}{4} x^2 + \frac{15}{4} x^3 + x^4$
can be written as
\begin{displaymath}
p_\eta(x) = \bigl(q_{\eta,1}(x)\bigr)^2 + \bigl(q_{\eta,2}(x)\bigr)^2
    + \ldots + \bigl(q_{\eta,n_\eta}(x)\bigr)^2
    = \bigl( \bar{q}(x)\bigr)\trp \bar{q}(x)
\end{displaymath}
for some $n_\eta \in \numN$ and polynomials $q_{\eta,i} \in \numR[x]$.

\begin{smremark}
Not every non-negative polynomial can be written as a sum of squares
\cite{Blekherman2010a}% X10103465 arxiv
. However, an ``approximation'' might be possible
\cite{Lasserre2007b,Lasserre2017a,Parrilo2003a}% MR2375528 0036-1445
. 
\end{smremark}

\begin{Example}\label{ex:nc.sdp}
Let $p_\eta(x) = \eta + \frac{13}{4} x^2 + \frac{15}{4} x^3 + x^4 \in \numR[x]$.
It can be written as
\begin{displaymath}
p_\eta(x) =
\begin{bmatrix}
1 & x & x^2
\end{bmatrix}\underbrace{%
\begin{bmatrix}
\tabstrut \eta & 0 & -\frac{\xi}{2} \\
\tabstrut 0 & \frac{13}{4} + \xi & \frac{15}{8} \\
\tabstrut -\frac{\xi}{2} & \frac{15}{8} & 1
\end{bmatrix}}_{=:Q(\xi,\eta)}
\begin{bmatrix}
1 \\ x \\ x^2
\end{bmatrix}.
\end{displaymath}
For $\eta=2$, and (approximately) $0.30775 < \xi < 1.39437$,
the matrix $Q(\xi,2)$ is \emph{positive definite}.
For $\eta=1$, the matrix $Q(0.5,1)$ is \emph{positive semidefinite},
and for $\eta<1$ it has at least one negative eigenvalue.
Thus the \emph{global minimum} of $p_2(x)$ is~$1$ (at $x=-2$)
which we can find by solving the semidefinite program
\begin{displaymath}
\min_{0 \ledef X \in \symmat_3 }
  \Bigl\{ \iprod{C}{X} :
    \iprod{A_1}{X} = b_1, \iprod{A_2}{X} = b_2
    \iprod{A_3}{X} = b_3, \iprod{A_4}{X} = b_4 \Bigr\}
\end{displaymath}
with
\begin{align*}
A_1 &= 
\begin{bmatrix}
. & 1 & . \\
1 & . & . \\
. & . & . 
\end{bmatrix},\quad
A_2 =
\begin{bmatrix}
. & . & 1 \\
. & 1 & . \\
1 & . & . 
\end{bmatrix},\quad
A_3 =
\begin{bmatrix}
. & . & . \\
. & . & 1 \\
. & 1 & . 
\end{bmatrix},\quad
A_4 = 
\begin{bmatrix}
. & . & . \\
. & . & . \\
. & . & 1 
\end{bmatrix},\\
b &=
\begin{bmatrix}
0 \\ \tabstrut\frac{13}{4} \\ \tabstrut\frac{15}{4} \\ 1
\end{bmatrix}
\quad\text{and}\quad
C = 
\begin{bmatrix}
1 & . & . \\
. & . & . \\
. & . & . 
\end{bmatrix}.
\end{align*}
We write $x = \tovector X$ and denote by
$A\in \numR^{m\times n^2}$ ($m=4$, $n=3$)
the matrix with rows
$(\tovector A_i)\trp$ to be able to write
$A \tovector(X) = b$ for the constraints.
\end{Example}

The following problem runs indicates that (in principal)
\emph{quadratic} convergence is possible,
\emph{independent} of the starting point.
The resulting path for the initial matrix with $\xi=0.75$ and $\eta=2$
is shown in Figure~\ref{fig:nc.sdp.205} (blue line). 
Table~\ref{tab:nc.sdp} shows the number of iterations compared
with some classical solvers and that finding an initial \emph{feasible}
matrix is an important issue.
For a basic discussion see Section~\ref{sec:nc.fim}.
It might be possible to develop problem specific methods
for initial iterates.
Notice in particular the right column in the outputs with
the (Frobenius) norm of $X\inv$
and the advantage of centered search directions
for the (approximate) solution.

\begin{verbatim}
octave:1> ncminex
octave:2> X = ncminsdp(A_sdp, b_sdp, C_sdp, Q(0.75,2));

NCMINSDP Search GCE      Version 0.99       December 2018       (C) KS
----------------------------------------------------------------------
Semidefinite Program: n=3, m=4, m_dof=1, m_min=1
      1-tau=5.00e-05, tol=2.00e-08, mu=0.000, maxit=20
cnt   alpha      tr(C*dX)      tr(C*X)        ||A*vec(X)-b||  ||X^-1||
  0   0.00e+00   -1.000e-01     2.000000000e+00     0.00e+00   2.5e+01
  1   9.88e+00   -9.876e-01     1.012402018e+00     2.40e-16   1.7e+03
  2   1.26e-02   -1.240e-02     1.000002493e+00     5.06e-16   8.4e+06
  3   2.01e-04   -2.493e-06     1.000000000e+00     5.06e-16   1.7e+11
  4   5.00e-05   -1.247e-10     1.000000000e+00     9.21e-16   3.9e+15
                                1.000000000000002e+00
octave:3> norm(X-X_sdp)
ans =    2.0251e-15
\end{verbatim}

\begin{verbatim}
octave:4> X = ncminsdp(A_sdp, b_sdp, C_sdp, Q(0.75,100));

NCMINSDP Search GCE      Version 0.99       December 2018       (C) KS
----------------------------------------------------------------------
Semidefinite Program: n=3, m=4, m_dof=1, m_min=1
      1-tau=5.00e-05, tol=2.00e-08, mu=0.000, maxit=20
cnt   alpha      tr(C*dX)      tr(C*X)        ||A*vec(X)-b||  ||X^-1||
  0   0.00e+00   -1.000e-01     1.000000000e+02     0.00e+00   1.0e+01
  1   9.88e+02   -9.878e+01     1.219011619e+00     3.60e-14   9.9e+01
  2   2.04e-03   -2.010e-01     1.017994134e+00     3.64e-14   1.2e+03
  3   8.95e-02   -1.799e-02     1.000001178e+00     3.64e-14   1.8e+07
  4   6.55e-05   -1.178e-06     1.000000000e+00     3.64e-14   3.6e+11
  5   5.00e-05   -5.892e-11     1.000000000e+00     3.70e-14   5.9e+15
                                1.000000000000041e+00
octave:5> norm(X-X_sdp)
ans =    4.8726e-14
\end{verbatim}

\begin{verbatim}
octave:6> X = ncminsdp(A_sdp, b_sdp, C_sdp)

NCMINSDP Search STD      Version 0.99       December 2018       (C) KS
----------------------------------------------------------------------
Semidefinite Program: n=3, m=4, m_dof=1, m_min=1
      1-tau=1.00e-06, tol=2.00e-08, mu=0.765, maxit=20
ini   alpha      min(eig(X))   tr(C*X)        ||A*vec(X)-b||  ||X^-1||
  1   5.76e-01    4.247e-01     1.000000000e+00     1.86e+00   2.0e+03
  2   9.13e-01    4.673e-02     1.205800795e+00     1.63e-01   4.4e+03
  3   1.98e+00    4.037e-02     1.986977186e+00     5.15e-16   3.3e+01
cnt   alpha      tr(C*dX)      tr(C*X)        ||A*vec(X)-b||  ||X^-1||
  0   0.00e+00   -1.000e-01     1.986977186e+00     5.15e-16   2.5e+01
  1   2.83e+00   -2.826e-01     1.704376627e+00     4.84e-16   9.4e+01
  2   6.72e+00   -6.720e-01     1.032384907e+00     6.38e-16   1.4e+03
  3   3.24e-01   -3.238e-02     1.000009826e+00     9.99e-16   4.3e+06
  4   9.83e-05   -9.825e-06     1.000000000e+00     4.58e-16   1.1e+11
  5   3.88e-09   -3.883e-10     1.000000000e+00     9.99e-16   1.1e+15
                                1.000000000000039e+00
X =
   1.0000e+00   7.0203e-18  -2.5000e-01
   7.0203e-18   3.7500e+00   1.8750e+00
  -2.5000e-01   1.8750e+00   1.0000e+00

octave:7> norm(X-X_sdp,'fro')
ans =    7.1906-08
\end{verbatim}

\begin{table}
\begin{center}
\begin{tabular}{llr}
Solver & Algorithm & Iterations \\\hline\tabstrut
SeDuMi 1.30 \cite{Sturm1999a}% MR1778433 1055-6788
  & 0 & 19 \\
  & 1 v-corr. & 12 \\
  & 2 $xz$-corr (PC) & 9 \\\hline\tabstrut
SDPA 7.3.9 \cite{SDPA2008}% ??? manual
  & PC & 13 \\\hline\tabstrut
\texttt{ncminsdp} (Section~\ref{sec:nc.imp.sdp})
  & standard & $3+5$ \\
  & symmetric & $3+5$ \\
  & algebraic centered search & $3+5$ \\
  & geometric centered search & $3+4$ \\\hline
\end{tabular}
\end{center}
\vspace{-2ex}
\caption{The number of iterations (that for finding a feasible
  starting matrix $+$ that for minimization)
  for the semidefinite program
  from Example~\ref{ex:nc.sdp}.
  ``PC'' stands for \emph{predictor-corrector}
  which implies more expensive steps.
  }
\label{tab:nc.sdp}
\end{table}

\section{Implementation NCMIN}\label{sec:nc.imp}

The following implementation in \textsc{Octave} \cite{OCTAVE2018}% ??? manual
\ is meant for educational purposes only without any warranties.
By changing the function \verb|mldivide| to \verb|lsqminnorm|
it can be used also in \textsc{Matlab} \cite{MATLAB2018}% ??? manual
. (Notice however that the output can be slightly different due
to rounding, etc.)
It consists of a simple standalone linear solver \verb|ncminlp|
with net less than 70 lines (see Section~\ref{sec:nc.imp.lp}
for the code and Section~\ref{sec:nc.lp} for a typical application)
and a semidefinite solver \verb|ncminsdp| (Section~\ref{sec:nc.imp.sdp})
building on \verb|ncmindof| (Section~\ref{sec:nc.imp.dof}) and
\verb|ncminini| (Section~\ref{sec:nc.imp.ini})
with net less than 210 lines together.
For a typical problem run of the latter see Section~\ref{sec:nc.sos}.

\begin{smremark}
Please keep in mind that the code is not streamlined to be
able to follow each step without jumping into a cascade of
functions, in particular to serve as a thread through the
theory. Resisting the temptation to further ``optimize''
it was not that easy. But already for sparse matrices
some algorithms for ``simple'' tasks like solving a
\emph{underdetermined} linear system
of equations change significantly and hand\-ling the differences
between \textsc{Octave} and \textsc{Matlab} would make the
code unnecessary complicated (to read). Those who are interested in
solving large scale problems might find 
\cite{Morikuni2015a}% MR3317782 0895-4798
\ helpful. For a general background in numerical linear
algebra we refer to 
\cite{Demmel1997a}% MR1463942 book
.
One typical ``inefficiency'' is that the inverse of $X$ is computed
twice, for example in lines~26 and~34 (Section~\ref{sec:nc.imp.ini})
or lines~80 and~88 (Section~\ref{sec:nc.imp.sdp}). But to really
attack that seriously one needs to think about approximating the
inverse (and the consequences for the solution) in the case of
sparse matrices (``fill in'', etc.)
\end{smremark}

\begin{smremark}
If something goes wrong or the result is not like
expected then one should check every single input argument
and try different parameters. Is the matrix $C$ symmetric?
Is the (primal) problem \emph{feasible}? Is it \emph{bounded}?
What is the rank (what are the singular values) of $A$?
Is the tolerance appropriate (especially for a badly conditioned
problem)? Is the steplength scaling $\tau$ too ``aggressive''?
\end{smremark}

\begin{smremark}
If everything looks nice, have a look on the
\emph{dual} problem. Is the duality gap zero?
What are the ranks (with respect to some tolerance) of
the (almost) optimal pair $X$ and $Z$?
Is the complementary primal-dual (almost) optimal
solution \emph{strict} complementary
\cite{Wei2010a}% MR2718693 0025-5610
?
Is the result (within some tolerance) \emph{independent}
of the initial starting value?
Yield the other search strategies the same result?
Is the result ``stable'' with respect to a change of
the parameters (although the number of iterations changes
significantly)? Is the maximal number of iterations
sufficiently large? Are there some ``problem specific''
(plausibility) checks?
\end{smremark}

\lstset{language=Octave,basicstyle=\ttfamily}
\lstset{numberblanklines=true, basewidth={0.5em,0.45em}}
\lstset{morecomment=[l][\color{ks-blue}]{\%}}
\lstset{numbers=left, numberstyle=\footnotesize, stepnumber=5, firstnumber=1}
\lstset{frame=lines,aboveskip=0.5ex,belowskip=1.5ex}

\subsection{NCMIN Examples}\label{sec:nc.imp.ex}

These tiny examples are for basic testing and detailled
investigations. The matrix \verb|C_ubd| in line~17 yields
an \emph{unbounded} problem.
The family of examples from 
\cite{Jansson2007a}% MR2377260 0036-1429
\ are useful to check the ``limits'' of a solver.
Although there are some similarities to the
example from Section~\ref{sec:nc.sos}, there
are no \emph{degrees of freedom} and hence
search directions cannot be constructed
(search types~3 and~4 in \verb|ncminsdp|).
\textsc{SeDuMi}
\cite{Sturm1999a}% MR1778433 1055-6788
\ provides an interface to the sparse format of SDPA
\cite{Yamashita2003a}% MR2019042 1055-6788
.

\begin{lstlisting}[firstline=7,label=lst:nc.min.ex,%
  title={\texttt{ncminex.m}}]
% LP: simple linear program [Lecture Notes LinOpt, TU Graz, 2005]
A_lp = [ -2, 1, 1, 0, 0;  -1, 2, 0, 1, 0;  1, 0, 0, 0, 1 ];
b_lp = [ 2, 7, 3 ]';
c_lp = [ -1, -2, 0, 0, 0 ]';
x_lp = [ 3, 5, 3, 0, 0 ]'; % Optimal solution

% SDP: semidefinite program (from SOS)
% Global minimization of p = 2 + 13/4*x^2 + 15/4*x^3 + x^4
A_sdp = [ 0, 1, 0, 1, 0, 0, 0, 0, 0; ...
          0, 0, 1, 0, 1, 0, 1, 0, 0; ...
          0, 0, 0, 0, 0, 1, 0, 1, 0; ...
          0, 0, 0, 0, 0, 0, 0, 0, 1 ];
b_sdp = [ 0, 3.25, 3.75, 1]';
C_sdp = diag([ 1, 0, 0]);
C_ubd = diag([ -1, 0, 0 ]); % Unbounded problem
X_sdp = [ 1, 0, -0.25; 0, 3.75, 1.875; -0.25, 1.875, 1 ];
Q = @(xi,eta) [ eta, 0, -xi/2;  ...
                0, b_sdp(2)+xi, b_sdp(3)/2; ...
                -xi/2, b_sdp(3)/2, b_sdp(4) ];

% Jansson--Chaykin--Keil [SIAM J. Numer. Anal. 2007]
A_jck = [ 0, -0.5, 0, -0.5, 0, 0, 0, 0, 0; ...
          1, 0, 0, 0, 0, 0, 0, 0, 0; ...
          0, 0, 1, 0, 0, 0, 1, 0, 0; ...
          0, 0, 0, 0, 0, 1, 0, 1, 0 ];
b_jck = @(epsilon) [1, epsilon, 0, 0]';
C_jck = @(delta) [ 0, 0.5, 0; 0.5, delta, 0; 0, 0, delta ];
X_jck = @(epsilon) [ epsilon, -1, 0; -1, 1/epsilon, 0;  0, 0, 0 ];

% Relaxation of a combinatorial problem
A_cmb = [ 1, 0, 0, 0, 0, 0, 0, 0, 0; ...
          0, -0.5, 0, -0.5, 1, 0, 0, 0, 0; ...
          0, 0, -0.5, 0, 0, 0, -0.5, 0, 1 ];
b_cmb = [1, 0, 0]';
C_cmb = diag([0, 1, -0.5]);
\end{lstlisting}

\subsection{NCMIN Linear Solver}\label{sec:nc.imp.lp}

For the theoretical setting see Section~\ref{sec:nc.lp}
with the main linear system \eqref{eqn:nc.lp.search}
corresponding to ``\verb|A_sys*x_sys=b_sys|'' in lines 60--63.
Search type~2 (line~7) is only for illustrating that a naive
approach ---not ensuring invertibility of the (diagonal) matrix
as a linearized condition---
can easily trigger wrong results.
There are several parts which need further investigation.

Firstly, finding a feasible initial vector with positive
entries. This is a simplified version of that presented
in Section~\ref{sec:nc.fim} and computes iteratively
a direction ``towards'' $Ax=b$ subject to the linearized
invertibility condition
\begin{displaymath}
\diag \tilde{x} \diag \Delta x + \diag x \diag \Delta \tilde x = 0
\end{displaymath}
until the maximal steplength $\alpha = \alpha_{\max}(x,\Delta x) > 1$
or the residuum $b-Ax$ is smaller than a predefined tolerance
(or the maximal number of iterations is reached).

And secondly, it would be interesting how a modified version
of search type~1 including a correction of the linearization
error (compare with classical \emph{predictor-corrector}
methods \cite{Mehrotra1992a}% MR1186163 1052-6234
) works. 

For a small linear (scheduling) program with $m=24$ and $n=80$
\texttt{ncminlp} needs
7~iterations for the initial iterate and
18~iterations for the minimization (with relative error
$10^{-8}$).
\textsc{SeDuMi} 1.30 
\cite{Sturm1999a}% MR1778433 1055-6788
\ needs 8~iterations for the same problem with standard settings
(\emph{predictor-corrector} method)
and 32~iterations with the ``not recommended''
algorithm~0.
Notice however that the latter is a highly developed program,
while \texttt{ncminlp} has a lot of potential for improvements:
Reasonable (relative) accuracy is reached quite fast (after a few
iterations) and it is possible to stop at any point
\emph{without} impact on the feasibility of the (approximate) solution.
Additionally, the linear systems of equations to solve
are smaller than for classical primal-dual interior-point methods \ldots

\begin{lstlisting}[firstline=7,label=lst:nc.min.lp,%
  title={\texttt{ncminlp.m}}]
% min_{x>=0} { c'*x : A*x=b }
% [ n*1 ] = ncminlp(m*n, m*1, n*1, n*1)
% rank(A) = m, 0 < tau < 1
function [x, info] = ncminlp(A, b, c, tau)

  % Parameters
  typ = 1; % (1) all directions, (2) naive approach (for testing only)
  txt_search = [ 'STD'; 'TST' ];
  maxit = 30;
  if nargin < 4
    tau = 1-1e-10;
  end
  tol = 2e-9;
  trCdX = -0.1;
  maxit_ini = 12;
  tau_ini = 0.99999;
  tol_res = 5e-15;
  
  % Initialization
  [m, n] = size(A);
  M_dir = rref(null(A)')';
  m_dir = size(M_dir,2);
  fmt_0 = '\nNCMINLP Search %s       Version %s       %s       (C) KS\n';
  txt_0 = '-----------------------------------';
  fmt_1 = '%s   %s      %s   tr(c*x)             ||A*x-b|| ||x.^-1||\n';
  fmt_2 = '%3i   %.2e   % .3e    % .9e     %.2e   %.1e\n';
  fprintf(1, fmt_0, txt_search(typ,:), '0.99', 'December 2018');
  fprintf(1, '%s%s\n', txt_0, txt_0);
  fprintf(1, 'Linear Program: n=%i, m=%i\n', n, m);
  fprintf(1, '      1-tau=%.2e, tol=%.2e, maxit=%i\n', 1-tau, tol, maxit);

  % Find initial feasible vector
  x = ones(n,1);
  fprintf(1, fmt_1, 'ini', 'alpha', 'min(x)     ');
  for cnt=1:maxit_ini
    x_inv = x.^-1;
    A_sys = [ A, zeros(m,n); ...
              diag(x_inv), diag(x) ];
    b_sys = [ (b-A*x); zeros(n,1) ];
    x_sys = mldivide(A_sys, b_sys);
    dx = x_sys(1:n);
    alpha = 1/max(-dx./x);
    x = x + min(1,tau_ini*alpha)*dx;
    eig_1 = min(x);
    nrm_res = norm(A*x-b);
    fprintf(1, fmt_2, cnt, alpha, eig_1, c'*x, nrm_res, norm(x.^-1));
    if (eig_1 > 0) && (nrm_res < tol_res)
      break
    end
  end
  info.iter = cnt;
  x_ini = x;

  fprintf(1, fmt_1, 'cnt', 'alpha', 'tr(c*dx)   ');
  fprintf(1, fmt_2, 0, 0, trCdX, c'*x, norm(A*x-b), norm(x.^-1));

  % Minimization
  for cnt=1:maxit
    x_inv = x.^-1;
    A_sys = [ diag(x_inv)*M_dir, diag(x); ...
              c'*M_dir, zeros(1,n) ];
    b_sys = [ zeros(n,1); trCdX ];
    x_sys = mldivide(A_sys, b_sys);
    if typ == 2
      x_sys = mldivide(c'*M_dir,trCdX);
    end
    dx = M_dir*x_sys(1:m_dir,1);
    alpha = 1/max(-dx./x);
    dx = tau*alpha*dx;
    x = x + dx;

    info.res = norm(A*x-b);
    fprintf(1, fmt_2, cnt, alpha, c'*dx, c'*x, info.res, norm(x.^-1));
    if abs(c'*dx) < tol
      break
    end
  end
  info.iter = info.iter + cnt;
  fprintf(1, '                               % .15e\n\n', c'*x);
end
\end{lstlisting}

\subsection{NCMIN Vector Space Decomposition}\label{sec:nc.imp.dof}

The basic theory is explained in Section~\ref{sec:nc.nla}.
Here, \verb|M_dof| corresponds to $M_\xi$,
\verb|M_min| to $M_\eta$
and \verb|M_nsy| to $M_\nu$.
The loop from line~21 to~27 takes care of distinguishing
directions between centering and minimization
and does some simple
``normalization'' of the vectors in $M_\xi$.
Notice the limitations of this (naive) approach
with respect to the use of \emph{dense} matrices
(in total) of size $n^2 \times n^2$.
Although this task is not directly related to ``preprocessing''
it is maybe worth to view it from a more general context
\cite{Cheung2013a,Malick2009a,Ramana1997b,Tuncel2012a}% MR3108430 incollection
.

\begin{lstlisting}[firstline=7,label=lst:nc.min.dof,%
  title={\texttt{ncmindof.m}}]
% [ n^2*m1, n^2*m2, n^2*m3 ] = ncmindof(m*n^2, m*1, n*n, 1*1)
function [M_dof, M_min, M_nsy] = ncmindof(A, b, C, tol)
  n = size(C,2);
  m_nsy = n*(n-1)/2;
  M_nsy = zeros(n^2,m_nsy);
  k = 0;
  for i=1:n
    for j=i+1:n
      xi = zeros(n,n);
      xi(i,j) = 1;
      xi(j,i) = -1;
      k = k+1;
      M_nsy(:,k) = reshape(xi,n^2,1);
    end
  end
  M_dof = null([A; M_nsy']);
  m_dof = rank(M_dof);
  M_dof = rref(M_dof')';
  c_min = reshape(C,n^2,1);
  flg_min = ones(m_dof,1);
  for k=1:m_dof
    lst_idx = find(abs(M_dof(:,k)) == max(abs(M_dof(:,k))));
    M_dof(:,k) = M_dof(:,k) / M_dof(lst_idx(1),k);
    if abs(c_min'*M_dof(:,k)) < tol
      flg_min(k) = 0;
    end
  end
  M_min = M_dof(:,find(flg_min==1));
  M_dof = M_dof(:,find(flg_min==0));
end
\end{lstlisting}

\subsection{NCMIN Initial Feasible Matrix}\label{sec:nc.imp.ini}

For details see Section~\ref{sec:nc.fim}.
For Example~\ref{ex:nc.sdp} one can subtract \verb|Q(0,0)| from
a matrix $X \in \symmat_3$ to read off $\xi$ in entry $(2,2)$
because $M_\xi\trp=[0,0,-\frac{1}{2},0,1,0,-\frac{1}{2},0,0]\trp$
and $\eta$ in entry $(1,1)$ because
$M_\eta\trp=[1,0,\ldots,0]\trp$.
In the following example we have $\xi \approx 0.767$ and $\eta\approx 1.987$:
\begin{verbatim}
octave:1> ncminex
octave:2> X_ini = ncminini(A_sdp, b_sdp, C_sdp, 5e-16);
octave:3> X_ini - Q(0,0)
ans =
   1.9870e+00  -1.1829e-16  -3.8329e-01
  -1.1829e-16   7.6659e-01   2.2204e-16
  -3.8329e-01   2.2204e-16  -1.1102e-16
\end{verbatim}

\begin{lstlisting}[firstline=7,label=lst:nc.min.ini,%
  title={\texttt{ncminini.m}}]
% [ n*n ] = ncminini(m*n^2, m*1, n*n, 1*1)
function [X, info] = ncminini(A, b, C, tol)
  
  % Parameters
  tau = 0.9995;
  mu = 0.0;
  tol_res = 5e-15;
  maxit = 15;

  % Initialization
  m = size(A,1);
  n = size(C,2);
  I_n = eye(n);
  [M_dof, M_min, M_nsy] = ncmindof(A, b, C, tol);
  M_dir = [M_dof, M_min, M_nsy];
  m_dir = size(M_dir,2);
  m_dof = size(M_dof,2);
  fmt_1 = '%s   %s      %s   tr(C*X)        ||A*vec(X)-b||  ||X^-1||\n';
  fmt_2 = '%3i   %.2e   % .3e    % .9e     %.2e   %.1e\n';
  fprintf(1, fmt_1, 'ini', 'alpha', 'min(eig(X))');

  X = I_n;
  x = reshape(X,n^2,1);

  for cnt=1:maxit
    X_inv = X^-1;
    scl_wrk = 10*norm(X_inv);
    A_sys = [ scl_wrk*A, zeros(m,n^2); ...
              kron(I_n,X_inv), kron(I_n,X) ];
    b_sys = [ scl_wrk*(b-A*x); zeros(n^2,1) ];
    x_sys = mldivide(A_sys, b_sys);
    dX = reshape(x_sys(1:n^2),n,n);
    dX = 0.5*(dX+dX');
    L_inv = inv(chol(X))';
    alpha = 1/max(eig(-L_inv*dX*L_inv'));
    X = X + min(1,tau*alpha)*dX;
    x = reshape(X,n^2,1);

    % Centering
    X_inv = X^-1;
    if m_dof > 0
      A_sys = [ kron(I_n,X_inv)*M_dof, kron(I_n,X) ];
      x_sys = mldivide(A_sys, (1-mu)*reshape(X_inv,n^2,1));
      dX = reshape(M_dof*x_sys(1:m_dof), n, n);
      X = X + 0.5*(dX+dX');
    else
      A_sys = [ kron(I_n,X_inv)*M_dir, kron(I_n,X) ];
      x_sys = mldivide(A_sys, (1-mu)*reshape(X_inv,n^2,1));
      dX = reshape(M_dir*x_sys(1:m_dir), n, n);
      L_inv = inv(chol(X))';
      dX = 0.5*(dX+dX');
      eig_n = max(eig(-L_inv*dX*L_inv'));
      beta = 1;
      if eig_n > 1
        beta = tau/eig_n;
      end
      X = X + beta*dX;
    end
    x = reshape(X,n^2,1);

    eig_1 = min(eig(X));
    nrm_res = norm(A*x-b);
    fprintf(1, fmt_2, cnt, alpha, eig_1, trace(C'*X), ...
      nrm_res, norm(X_inv,'fro'));
    if (eig_1 > 0) && (nrm_res < tol_res)
      break
    end
  end
  info.iter = cnt;
end
\end{lstlisting}

\subsection{NCMIN Semidefinite Solver}\label{sec:nc.imp.sdp}

For details see Section~\ref{sec:nc.fip}.
For a ``typical'' path for search type~1 (all directions)
see Figure~\ref{fig:nc.sdp.201},
type~2 (symmetric directions) Figure~\ref{fig:nc.sdp.203},
type~3 (algebraic centered) Figure~\ref{fig:nc.sdp.204} and
type~4 (geometric centered) Figure~\ref{fig:nc.sdp.205}.
In any case it is a \emph{greedy} method, that is,
after computing a ``feasible'' search direction,
it takes (almost) the maximal steplength.
The main ingredient for making that work is
\emph{centering}, or in other words, ``staying away from the boundary (of
singular matrices)''. See Section~\ref{sec:nc.nla},
in particular Example~\ref{ex:nc.geocen} and~\ref{ex:nc.algcen}.

\begin{lstlisting}[firstline=7,label=lst:nc.min.sdp,%
  title={\texttt{ncminsdp.m}}]
% min_{X>=0} { tr(C'*X) : A*vec(X)=b }
% [ n*n, n^2*k ] = ncminsdp(m*n^2, m*1, n*n, n*n, 1*1)
% rank(A) = m, 0 <= mu <= 1, X_ini > 0
function [X, info, X_path] = ncminsdp(A, b, C, X_ini, mu)

  % Parameters
  typ = 1; % (1) all directions, (2) symmetric directions,
           % constructed by algebraic (3) or geometric (4) centering
  txt_lst = [ 'STD'; 'SYM'; 'ACE'; 'GCE' ];
  mu_lst = [ 0.765, 0.235, 0.215, 0.0 ];
  tau_lst = [ 0.999999, 0.99995, 0.99999, 0.99995 ];
  maxit = 20;         % Maximal number of iterations
  tol = 2e-8;         % Tolerance
  tau = tau_lst(typ); % Steplength scaling tau < 1
  if nargin < 5
    mu = mu_lst(typ); % Centering parameter 0 <= mu <= 1
  end
  trCdX = -0.1;

  % Initialization
  m = size(A,1);
  n = size(C,2);
  I_n = eye(n);
  C_mtx = reshape(C', n^2, 1)';
  [M_dof, M_min, M_nsy] = ncmindof(A, b, C, 5e-16);
  m_dof = size(M_dof,2);
  M_dir = [M_dof, M_min];
  if typ == 1
    M_dir = [M_dir, M_nsy];
  end
  m_dir = size(M_dir,2);
  fmt_0 = '\nNCMINSDP Search %s      Version %s       %s       (C) KS\n';
  txt_0 = '-----------------------------------';
  fmt_1 = '%s   %s      %s   tr(C*X)        ||A*vec(X)-b||  ||X^-1||\n';
  fmt_2 = '%3i   %.2e   % .3e    % .9e     %.2e   %.1e\n';
  fprintf(1, fmt_0, txt_lst(typ,:), '0.99', 'December 2018');
  fprintf(1, '%s%s\n', txt_0, txt_0);
  fprintf(1, 'Semidefinite Program: n=%i, m=%i, m_dof=%i, m_min=%i\n', ...
    n, m, size(M_dof,2), size(M_min,2));
  fprintf(1, '      1-tau=%.2e, tol=%.2e, mu=%.3f, maxit=%i\n', ...
    1-tau, tol, mu, maxit);
  if nargin < 4
    [X_ini, info_ini] = ncminini(A, b, C, 5e-16);
  else
    info_ini.iter = 0;
  end
  X = X_ini;
  x = reshape(X,n^2,1);
  X_path = x;
  fprintf(1, fmt_1, 'cnt', 'alpha', 'tr(C*dX)   ');
  fprintf(1, fmt_2, 0, 0, trCdX, trace(C'*X), norm(A*x-b), norm(X^-1,'fro'));
  
  % Centering
  if typ >= 3
    X_inv = X^-1;
    A_sys = [ kron(I_n,X_inv)*M_dof, kron(I_n,X) ];
    x_sys = mldivide(A_sys, (1-mu)*reshape(X_inv,n^2,1));
    dX = reshape(M_dof*x_sys(1:m_dof), n, n);
    X_cen = X + 0.5*(dX+dX');
    if typ == 4
      L_inv = inv(chol(X_cen))';
      X_new = zeros(n,n);
      for i=1:m_dof
        dX = reshape(M_dof(:,i),n,n);
        alpha_1 = 1/max(eig(-L_inv*dX*L_inv'));
        alpha_2 = 1/max(eig(L_inv*dX*L_inv'));
        X_new = X_new + 0.5*(alpha_1-alpha_2)*dX;
      end
      X_cen = X_cen + (1-mu)*X_new;
    end
  end

  % Minimization
  for cnt=1:maxit
    X_old = X;
    if (typ >= 3) && (cnt > 1)
      dX = X - X_cen;
      X_cen = X;
    else
      X_inv = X^-1;
      A_sys = [ kron(I_n,X_inv)*M_dir, kron(I_n,X); ...
                C_mtx*M_dir, zeros(1,n^2) ];
      b_sys = [ zeros(n^2,1); trCdX ];
      x_sys = mldivide(A_sys, b_sys);
      dX = reshape(M_dir*x_sys(1:m_dir), n, n);
      dX = 0.5*(dX+dX');
    end
    L_inv = inv(chol(X))';
    eig_n = max(eig(-L_inv*dX*L_inv'));
    if abs(eig_n) < 1e-20
      fprintf(1, 'Warning: Problem possibly unbounded.\n');
      alpha = 1e4;
    else
      alpha = 1/eig_n;
    end
    X = X + tau*alpha*dX;
    x = reshape(X,n^2,1);
    X_path = [X_path, x];
    
    % Centering
    if (m_dof > 0) && (abs(trace(C'*tau*alpha*dX)) > 1e-3*tol)
      X_inv = X^-1;
      A_sys = [ kron(I_n,X_inv)*M_dof, kron(I_n,X) ];
      x_sys = mldivide(A_sys, (1-mu)*reshape(X_inv,n^2,1));
      dX = reshape(M_dof*x_sys(1:m_dof), n, n);
      eig_tmp = eig(X)';
      X = X + 0.5*(dX+dX');
      if typ == 4
        L_inv = inv(chol(X))';
        X_new = zeros(n,n);
        for i=1:m_dof
          dX = reshape(M_dof(:,i),n,n);
          alpha_1 = 1/max(eig(-L_inv*dX*L_inv'));
          alpha_2 = 1/max(eig(L_inv*dX*L_inv'));
          X_new = X_new + 0.5*(alpha_1-alpha_2)*dX;
        end
        X = X + (1-mu)*X_new;
      end
      x = reshape(X,n^2,1);
      X_path = [X_path, x];
    end

    tr_1 = trace(C'*(X-X_old));
    info.res = norm(A*x-b);
    nrm_inv = norm(X^-1, 'fro');
    fprintf(1, fmt_2, cnt, alpha, tr_1, trace(C'*X), info.res, nrm_inv);
    if abs(tr_1) < tol
      break
    end
  end
  info.iter = info_ini.iter + cnt;
  fprintf(1, '                               % .15e\n\n', trace(C'*X));
end
\end{lstlisting}

\section{(Not) Only for Students}\label{sec:nc.stud}

Starting with semidefinite programming \emph{is} difficult,
even with a ``linear'' background. Maybe because there is
no end in sight and sooner than later one is in the middle
of current research or needs deep results which are hard
to understand even roughly within a few weeks (or months).
Just to get an impression:

``Semidefinite programming (SDP) is probably the most important
new development in optimization in the last two decades. (\ldots)
Notice that semidefinite programming is a far reaching extension
of linear programming (LP) \ldots'' 
\cite{Vinnikov2012b}% MR2962792 incollection
.

Its importance cannot be overestimated due to the many
practical applications, in particular in engineering
\cite[Section~4.3]{Nemirovski2007a}% MR2334199 incollection
\ \emph{and} the availability of highly developed
solvers. But semidefinite programming is also fascinating
from a purely mathematical point of view because of all
the manifold connections to (at a first glance) different
mathematical areas and the rich ``explorable'' structure.
The initial hurdle of many definitions and different notations
is admittedly high. Getting confident in (applied)
non-commutative algebra takes some time \ldots

So, \emph{where} should one start? With \emph{small}
problems that can be ``touched'' and from a perspective
one is familiar with, no matter if it is of geometrical,
analytical, algebraic or applied (in the sense of programming)
nature. Almost all here originated from a ``tiny'' semidefinite
problem with $3 \times 3$ matrices. And if there were not
Halmos' recommendation to \emph{stop}
\cite[Section~19]{Halmos1970b}% MR0277319 0013-8584
, this paper would have never been finished (and therefore not
exist).

\medskip
One could start with the simple \emph{linear matrix pencil}
$A = A(x,y) \in \numR^{2 \times 2}$
from \cite[Example~1.1.1]{Netzer2011a}% TH106 thesis
\ and try to figure out, which ``object'' it describes
with respect to the \emph{real} variables/parameters
$x$ and $y$ such that $A$ is \emph{positive (semi-)definite}:
\begin{displaymath}
A = A(x,y) = 
\begin{bmatrix}
1+x & y \\
y & 1-x
\end{bmatrix}
=
\underbrace{%
\begin{bmatrix}
1 & . \\
. & 1
\end{bmatrix}}_{=:A_0}
\otimes 1 +
\underbrace{%
\begin{bmatrix}
1 & . \\
. & -1
\end{bmatrix}}_{=:A_x}
\otimes x +
\underbrace{%
\begin{bmatrix}
. & 1 \\
1 & .
\end{bmatrix}}_{=:A_y}
\otimes y.
\end{displaymath}
Here the zeros are replaced by (lower) dots to emphasize the
structure and the \emph{tensor product} ``$\otimes$'' just
means that we ``plug in'' the variables (into the respective
matrix). Depending on the context/area, matrix pencils are
written in many different ways, for example
$A = A(x,y) = A_0 + x A_x + y A_y$.
The ``constant'' coefficient matrix $A_0$ usually plays
a special role (and is often just the identity matrix).

How can we find a matrix $X_0$
``inside'' this \emph{convex} object?
And given $X_0$, how can we find the ``center''
algorithmically by using \eqref{eqn:nc.nla.steplength}
respectively lines~61ff or~86ff in \verb|ncminsdp|
(Section~\ref{sec:nc.imp.sdp})?
Which other objects can be represented in such a way
(by $2 \times 2$ matrices)?
How many variables are needed at most for $3 \times 3$ matrices?
(Recall that we are talking about \emph{symmetric} matrices only.)

Questions about the representation of convex sets by
linear matrix pencils (or ``linear matrix inequalities'', LMI's)
are very difficult in general. But it is still worth to get
an idea by having a look on the many illustrations in
\cite{Netzer2011a}% TH106 thesis
.

\medskip
Another possibility is to use the power of semidefinite
programming for concrete problems. One application
is discussed in Section~\ref{sec:nc.sos}.
Another is the \emph{relaxation} of
hard combinatorial problems
\cite[Section~4.3.2]{Nemirovski2007a}% MR2334199 incollection
. For an illustration we restrict ourselves to 0-1 variables
and a quadratic functional $f : \numR^2 \to \numR$.
Let $C \in \numR^{3 \times 3}$,
\begin{displaymath}
A_0 =
\begin{bmatrix}
1 & . & . \\
. & . & . \\
. & . & .
\end{bmatrix},
\quad 
A_1 =
\begin{bmatrix}
. & -\frac{1}{2} & . \\
-\frac{1}{2} & 1 & . \\
. & . & .
\end{bmatrix}
\quad\text{and}\quad
A_2 = \begin{bmatrix}
. & . & -\frac{1}{2} \\
. & . & .  \\
-\frac{1}{2} & . & 1
\end{bmatrix}
\end{displaymath}
and define $\bar{x} = [ 1, x_1, x_2 ]\trp$
and $\bar{X} = \bar{x}\trp \bar{x}$.
Now the optimum from the (non-convex) combinatorial problem
can be bounded from below by a semidefinite program:
\begin{align}\label{eqn:nc.comb}
\min_{x \in \numR^2} \Bigl\{ f(x) : x_i \in \{ 0,1 \} \Bigr\} & \nonumber\\
\ge \min_{0\ledef X \in \symmat_3}  \Bigl\{  \iprod{C}{X} \Bigr.
  &: \Bigl. \iprod{A_0}{X} = 1, \iprod{A_1}{X} = 0, \iprod{A_2}{X} = 0 \Bigr\}.
\end{align}
Without the application of the trace the constraints
reads (typically)
\begin{displaymath}
A_1 X = 
\begin{bmatrix}
. & -\frac{1}{2} & . \\
-\frac{1}{2} & 1 & . \\
. & . & .
\end{bmatrix}
\begin{bmatrix}
1 & x_1 & x_2 \\
x_1 & x_1^2 & x_1 x_2 \\
x_2 & x_1 x_2 & x_2^2
\end{bmatrix}
=
\begin{bmatrix}
-\frac{x_1}{2} & * & * \\
* & -\frac{x_1}{2} + x_1^2 & * \\
* & * & 0
\end{bmatrix}.
\end{displaymath}
Some trivial but \emph{crucial} questions:
What is the rank of $\bar{X}$? What that
of $X \gnedef 0$ in the SDP~\eqref{eqn:nc.comb}?
And which linear algebraic concept(s) are useful
to measure the (approximate) ``numerical'' rank
of some $X \gnedef 0$?

The following works only for \emph{symmetric}
search directions (\verb|typ=2| in line~7)
and \emph{deactivated} centering (by uncommenting
line~106 in Section~\ref{sec:nc.imp.sdp} and
line~45 in Section~\ref{sec:nc.imp.ini}).
There might be a simple explanation \ldots

\begin{verbatim}
octave:1> ncminex; X = ncminsdp(A_cmb, b_cmb, C_cmb, X_ini)

NCMINSDP Search SYM      Version 0.99       December 2018       (C) KS
----------------------------------------------------------------------
Semidefinite Program: n=3, m=3, m_dof=1, m_min=2
      1-tau=5.00e-05, tol=2.00e-08, mu=0.235, maxit=20
ini   alpha      min(eig(X))   tr(C*X)        ||A*vec(X)-b||  ||X^-1||
  1   1.10e+00    8.932e-02     1.666666667e-01     5.55e-16   1.2e+01
cnt   alpha      tr(C*dX)      tr(C*X)        ||A*vec(X)-b||  ||X^-1||
  0   0.00e+00   -1.000e-01     1.666666667e-01     5.55e-16   1.2e+01
  1   4.95e+00   -4.950e-01    -3.283080272e-01     5.90e-16   6.0e+04
  2   1.70e+00   -1.702e-01    -4.985070477e-01     6.05e-16   1.8e+05
  3   1.26e-04   -1.258e-05    -4.985196251e-01     6.05e-16   3.5e+09
  4   7.03e+10   -6.267e-10    -4.985196257e-01     6.05e-16   7.1e+13
                               -4.985196257109498e-01
X =
   1.0000e+00   9.8043e-04   9.9900e-01
   9.8043e-04   9.8043e-04   1.9686e-03
   9.9900e-01   1.9686e-03   9.9900e-01
\end{verbatim}

\section{Epilogue}\label{sec:nc.epi}

Or \emph{pro}logue? Is there enough here that makes it worth
to look back? At the end everything grew out naturally 
from the trivial observation that the matrix multiplication is non-commutative
(in general). But what is surprisingly simple from a latter perspective is
not at all obvious in the beginning. And what we know from an old
proverb, namely that the search for a needle in a haystack
does not get easier when one increases the stack, seems to
hold true in particular for a ``matrix valued'' needle
(pointing towards the minimum). So one challenging question
remains: How can we find ``good''
search directions ---as solutions to underdetermined
\emph{linear} systems of equations--- along a
\emph{feasible} interior path?

Here we indicated only that it \emph{might} be possible
by combining some rather simple techniques. But even if it turns
out that these are not usable in general it could stimulate
new ideas and might explain some ``strange'' behaviour
\cite{Waki2012a}% MR3001030 0926-6003
\ or ``hard'' problem 
\cite{Wei2010a}% MR2718693 0025-5610
.
In any case one should keep Todd's words
from the abstract of \cite{Todd2001a}% MR2009698 0962-4929
\ in mind:
``The most effective computational methods are not always
provable efficient in theory, and \emph{vice versa}.''
Very interesting to read is also 
\cite[Section~1.1]{Forsgren2002a}% MR1980444 0036-1445
\ from Forsgren, Gill and Wright about the roots of
non-linear programming:
``Furthermore, a simplex-centric world view had the effect that
even `new' techniques mimicked the motivation of the simplex method
by always staying on a subset of exactly satisfied constraints.''
With other words: When is it necessary to \emph{leave} a
\emph{feasible} path?

\medskip
``Pedagogical and philosophical issues remain about the
best way to motivate interior-point methods ---perturbing optimality
conditions? minimizing a barrier function?--- and the multiplicity
of viewpoints continues to create new insights and new algorithms.''
\cite[Section~1.1]{Forsgren2002a}% MR1980444 0036-1445

\medskip
There is quite a lot we could not even touch here.
So we just mention a non-complete list
of papers (of quite different flavour)
which might help to find further literature.
For an overview and/or introduction
one could start with
\cite{Nemirovski2007a,Todd2001a,Vandenberghe1996a,Wright2005a}% MR2009698 0962-4929
.
Although one should use the latest versions (and the
current documentation) of the implementations, the
underlying publications
\cite{Toh1999a,Sturm1999a,Yamashita2003a}% MR1778429 1055-6788
\ are in some sense timeless.
More on interior-point methods, search directions,
predictor-corrector methods, etc.~can be found in
\cite{Alizadeh1998a,Helmberg1996a,Nemirovski2008a,Potra2000a,Todd1998a,Toh2002a}% MR1803305 0377-0427
.
Around verification, exact solutions, high precision, etc.~there are
\cite{Henrion2016a,Jansson2007a,Waki2012a}% MR3574590 1052-6234
\ and
\cite[Section~14]{Rump2010a}% MR2652784 0962-4929
.
Around strong duality, regularization, preprocessing,
hard and bad semidefinite programs, etc.~one should have a look in 
\cite{Cheung2013a,Malick2009a,Pataki2017a,Ramana1997a,Tuncel2012a,Wei2010a}% MR2988175 0926-6003
.

On one hand, semidefinite programming is a special
case of \emph{conic} programming
\cite{Belloni2009a,Boyd2004a,Forsgren2002a,Freund1999a,Jarre1992a,Jarre1995a,Nesterov1998a,Nesterov2012a}% MR1175483 0095-4616
.
On the other, it
is very useful for relaxations of hard combinatorial
problems. So there is a lot which goes far beyond
the basics we presented here:
\cite{Alizadeh1995a,Boyd1994a,Henrion2003a,Marshall2003a,Parrilo2003a}% MR1994866 0008-4395
, \cite[Section~4]{Nemirovski2007a}% MR2334199 incollection
.
Since semidefinite programming is tight together
with numerics (least squares approximation, iterative
linear solvers, finite precision arithmetics, etc.)~one
should mention at least
\cite{Cui2016b,Morikuni2015a,Waki2012a}% MR3317782 0895-4798
. But semidefinite programming is also interwoven
with \emph{semialgebraic geometry} 
\cite{Helton2010a,Netzer2011a,Helton2012a}% TH106 thesis
.

\medskip
At a first glance a topic like
``non-commutative sum of squares''
\cite{Helton2002a,Cafuta2011a}% MR1933721 0003-486X
\ seems to be far away. But indeed it was
one of the initial questions around
\emph{non-commutative convexity} 
\cite{Helton2006a}% MR2259894 0022-1236
\ which spans the bridge to my
``non-commutative'' algebraic research.

\medskip
There is a lot which is not discussed here (in detail).
One topic is numerics (in general) and sparse matrices
(in particular) in combination with iterative linear
solvers. One starting point for literature is
\cite{Morikuni2015a}% MR3317782 0895-4798
. Another is that around estimating good (maybe problem
specific) parameters. Since there are not that many in
the presented methods it should not be that difficult
to change them dynamically. One advantage is that there
is no need for a subtle control for $\mu \to 0^+$.
One disadvantage is that there is no continuous path
(and no analysis) anymore, so it might be difficult to
estimate a rate of convergence.

It is clear that the different methods can be combined,
in particular with ``classical'' methods. Additionally,
information from the dual problem can be used for
verification
\cite[Section~14]{Rump2010a}% MR2652784 0962-4929
. That this should not be underestimated is shown
in particular in 
\cite{Waki2012a}% MR3001030 0926-6003
.

The list of (not just small) improvements is long:
enable the use of a relative error and/or numerical limits 
as stopping criteria, export the initial data or some
intermediate results in \textsc{SDPA}'s 
\cite{SDPA2008}% ??? manual
\ sparse format, etc.
For the issues with respect to benchmarking
one could start with \cite{Mittelmann2003a}% MR1976487 0025-5610
, for large scale programs with
\cite{Fujisawa2007a}% MR2388512 0453-4514
.

In addition, thorough testing is necessary,
in particular with the problems from SDPLIB
\cite{Borchers1999a}% MR1778436 1055-6788
\ and the library \cite{deKlerk2009a}% MR2569727 1055-6788
. However, even if some basic testing could be done
automatically, a deeper analysis of the reasons
where something goes wrong (with the presented approach
or some classical solver) will definitely help to
get further insight in semidefinite programming.
In this context the survey
\cite{deKlerk2010a}% MR2573325 0377-2217
\ could be interesting.

\section*{Acknowledgement}

I am very grateful for Roland Speicher's invitation to
Saarbrücken and the opportunity to give a talk on this
topic and continue in particular ``non-commutative''
algebraic discussions in October 2018.
I thank Gabor Pataki and Michael J.~Todd for hints on the
literature and Birgit Janko for very valuable feedback on
several drafts of this work.

This work has been supported by research subsidies granted by
the government of Upper Austria
(research project ``Methodenentwicklung für Energieflussoptimierung'').

\addcontentsline{toc}{section}{Bibliography}
\bibliographystyle{alpha}
\bibliography{doku}

\end{document}